\documentclass[12pt,a4paper]{article}

\usepackage{amsmath}
\usepackage{amssymb}
\usepackage{amscd}

  \newcommand{\bs}[1]{\boldsymbol{#1}}

  \newcommand{\TT}{\mathbb T}
  \newcommand{\NN}{\mathbb N}
  \newcommand{\RR}{\mathbb R}
  \newcommand{\CC}{\mathbb C}
  \newcommand{\ZZ}{\mathbb Z}
  \newcommand{\QQ}{\mathbb Q}

  \newcommand{\oplam}{\mbox{\large $\curlywedge$}}
  \newcommand{\qed}{\hfill $\square$}
  \newcommand{\es}{\varnothing}
  \newcommand{\cA}{{\mathcal A}}
  \newcommand{\cK}{{\mathcal K}}
  \newcommand{\cM}{{\mathcal M}}
  \newcommand{\gD}{\varDelta}
  \newcommand{\gL}{\varLambda}
  \newcommand{\gG}{\varGamma}
  \newcommand{\ve}{\varepsilon}

  \newtheorem{theorem}{Theorem}
  \newtheorem{lemma}{Lemma}
  \newtheorem{prop}{Proposition}
  \newtheorem{coro}{Corollary}
  \newtheorem{fact}{Fact}
  \newtheorem{axiom}{Axiom}

\begin{document}

\centerline{\large \bf Weighted Dirac combs with pure point diffraction}

\vspace{15mm}

\centerline{{\sc Michael Baake$^{\,\rm a}$}
         and {\sc Robert V.\ Moody}$^{\,\rm b}$}
\vspace{10mm}

{\small
\hspace*{4em}
a:  Institut f\"ur Mathematik, Universit\"at Greifswald, \\
\hspace*{4em}
     \hspace*{2.5em} Jahnstr.\ 15a, 17487 Greifswald, 
      Germany\footnote{Address from 1.2.2004: Fakult\"at f\"ur
      Mathematik, Univ.\ Bielefeld, Postfach 100131, 33501
      Bielefeld, Germany}

\smallskip
\hspace*{4em}
b:  Dept. of Mathematical Sciences, University of Alberta, \\
\hspace*{4em}
     \hspace*{2.5em} Edmonton, Alberta T6G 2G1, Canada  }

\vspace{20mm}
\begin{abstract}
A class of translation bounded complex measures,  which have the
form of weighted Dirac combs, on locally compact Abelian groups
is investigated. Given such a Dirac comb, we are interested in its
diffraction spectrum which emerges as the Fourier transform of the
autocorrelation measure.
We present a sufficient set of conditions to ensure that the
diffraction measure is a pure point measure. Simultaneously,
we establish a natural link to the cut and project
formalism and to the theory of almost periodic measures.
Our conditions are general enough to cover the known theory of
model sets, but also to include examples such as the visible
lattice points.
\end{abstract}

\vspace{20mm}
\hphantom{Key Words:}
            almost periodic and pure point measures

\smallskip
\noindent
2000 MCS: $\;$42A25, 52C23 (primary); 43A05, 22B05, 28A33  (secondary)

\clearpage
\section{Introduction} \label{sec0}

The diffraction properties of quasicrystals, both physical
and mathematical, are among their most striking features.
The customary ways of explaining diffractivity are principally
of two kinds. The first is to view the quasicrystal as
orginating from a cut through a higher dimensional lattice,
in which case diffraction is explained as the vestiges of
the diffractive properties of the lattice that survive after
projection, see \cite{Hof,Duneau} and references therein.
The second is to look at the dynamical system constructed from
the hull of the quasicrystal under translation and to connect
the diffraction with the corresponding dynamical spectrum
\cite{Boris,Martin2,LMS-1}. In the latter case, an important r\^ole is
played by the repetitivity and unique ergodicity of the basic point
set (or tiling) that makes up the quasicrystal in question.
Despite serious efforts \cite{Sol,LM,LMS-2}, no complete answer is
known on how these two approaches fit together.

In this paper, we follow a different approach that accounts for
the diffractivity through three basic properties of the underlying
point set (or distribution of density). It is mathematically quite
direct, has the advantage of clarifying where internal spaces and
the higher dimensional lattices come from, and is
intrinsically insensitive to density $0$ changes.
Furthermore, the result will be considerably more general
than previous ones while the conditions can actually be verified
explicitly in many relevant examples.

Before we start our analysis, let us briefly comment on the notion of
diffraction used in this paper. Mathematically, as will become clear
in a moment, we are dealing with certain spectral properties of
measures, and the paper is rather self-contained in this respect.
Physically, our question is related to the X-ray diffraction image
of a solid, see \cite{Cowley} for background information. However,
we are completely ignoring the question whether the measure under
investigation is, in any sense, a realistic physical structure, such
as the ground state of some given atomic interaction, compare
\cite{Enter,BH} and references cited there. We simply assume that
it is given, and start from there. In particular, we do {\em not\/}
assume that the measures we will be considering are typical in any
sense. This would require us to couple our approach to ergodic theory
and to Gibbs measures, something which we defer for future work.

Also, strictly speaking, our results would only apply to a static
situation such as that of a solid at zero temperature (or to any
situation which is sufficiently well approximated by this point of
view). This poses no problem because, in a second step, one can
always extend the results rigorously to cover the dynamic situation
of diffraction at high temperatures, see \cite{Hof2} for details.
In view of this, it is perfectly adequate to withdraw from the
physical reality and to study the spectral properties of certain
measures in their own right, without losing the possible applications
to physics and crystallography.

\smallskip
The basic object of our study is then a translation bounded
(complex) measure $\omega$ on a locally compact Abelian group
$G$. In the standard case of a mathematical quasicrystal,
represented by a point set $\gL$,
$G$ would be the physical space $\RR^n$ and $\omega$ would
be the point measure $\sum_{x\in \gL} \delta_x$ which
is the Dirac comb supported on $\gL$. The basic assumptions
are three:
\begin{enumerate}
\item that the averaged autocorrelation $\gamma^{}_{\omega}$ of $\omega$
       should exist,
\item that the support of $\gamma^{}_{\omega}$ should be a uniformly
       discrete set $\gD$, and
\item that the $\ve$-almost periods of $\gamma^{}_{\omega}$
       should be relatively dense.
\end{enumerate}

Under these assumptions, we prove that the Fourier transform
$\hat{\gamma}^{}_{\omega}$ of the autocorrelation measure
$\gamma^{}_{\omega}$
is a pure point measure on $\hat{G}$, the dual group of $G$.
This is precisely the mathematical meaning of saying that
$\omega$ is pure point diffractive. We also show that all model
sets (no matter what internal spaces are involved) based on
windows with boundaries of (Haar) measure $0$ satisfy these
hypotheses, thus reproving the various results of this type
(cf.\ \cite{Martin2} and references given there) in a new and simpler
way. Here, the crucial link is made through Weyl's theorem on uniform
distribution. Moreover, there are interesting situations which are not
model sets, notably the visible points of a lattice (which do
{\em not\/} form a Delone set because they fail to be relatively dense),
that are also covered by our theorem. This
is a simplification of the previous approach in \cite{BMP}.

The constructions that are used in the proof of the main theorem
are as interesting to us as the result itself.
We introduce first a new group, namely the subgroup $L$ of $G$ generated
by the support of the autocorrelation measure $\gamma^{}_{\omega}$, and
second a topology (in fact, a uniformity) on $L$ that is derived from
the autocorrelation itself.
The latter is, in general, totally different from the one that $L$ gets
induced on it from $G$. The completion of $L$ with respect
to this autocorrelation topology is a new locally compact
Abelian group $H$. It is this group $H$ that is the `internal
group' of the system. Although $L$ is, in general, neither
discrete nor closed in $G$, it can be mapped into
$G \times H$ by essentially embedding it diagonally, whereupon
it appears as a {\em lattice}, denoted by $\tilde L$. Thus the
internal space
and the lattice, which seem normally to be without physical
meaning, appear naturally as a reflection of the $\ve$-almost
periodicity of the basic measure $\omega$, as defined through
its autocorrelation.

Furthermore, $G$ itself can be given a new topology which
combines its usual topology with the autocorrelation
topology of $L$ which lies inside it. The completion of
$G$ with respect to this new topology is a {\em compact} group
and is, in fact, nothing more than $(G \times H) / \tilde L$.
This gives a simple explanation of the appearance of this
compact group in the cut and project formalism which, though
it has been used a number of times in the so-called
`torus parameterization' \cite{BHP,HRB,Martin2}, did not seem
to have any {\em intrinsic\/} meaning before.

The appearance of $\ve$-almost periodicity is strongly
reminiscent of the theory of almost periodic functions of
Harold Bohr, though now it appears in the generalized form
of measures. The concept of almost periodicity of measures
in the theory of quasicrystals has appeared before \cite{Boris}.
In fact, during the work for this article, we became aware that the Bohr
theory has already been completely generalized to measures in a
fundamental but rather neglected paper of Gil de Lamadrid and Argabright
\cite{GdeL}. In some sense, this work already partly contains the
main result of this paper as will become clear in Section 
\ref{extensions}.
However, our method to arrive there, for a restricted but relevant
class of measures, is constructive and far more direct. Furthermore,
it has the virtue, as we have mentioned,
of making many of the common features in quasicrystal theory
appear by themselves in a totally natural way.

The setting of this paper is that of locally compact Abelian (LCA) 
groups
which are also $\sigma$-compact, and that of possibly unbounded, but
translation bounded, complex measures. The importance of LCA groups
for internal spaces of model sets has already been demonstrated,
see \cite{BMS} for concrete examples. The need for $\sigma$-compact 
LCA groups for the setting of quasicrystals themselves 
is not so obvious, but has been used before in \cite{Martin1}.
We have adopted it here because it does not pose extra complications
and it is a natural setting for the theorem and its proof. Readers who
prefer a more concrete base can replace $G$ by $\RR^n$ and the measure
$\omega$ by a discretely supported Dirac comb throughout.

\section{General setup and class of measures} \label{sec1}

Let $G$ be a $\sigma$-compact, locally compact Abelian (LCA) group, which
we always understand to include the Hausdorff property. The group action
will be written additively. The group $G$ is equipped
with a fixed Haar measure $\theta$ which is unique up to a scalar
multiple. If necessary, we will write $\theta^{}_G$ to link the Haar
measure to its group. Later on, we will need the {\em dual group\/},
denoted by $\hat{G}$, and we will then assume that it is equipped with
a matching Haar measure $\theta^{}_{\hat{G}}$ such that the generalized
Plancherel formula holds, see \cite[Thm.\ 2.5]{BF}. In this paper, the
word `measure' will always mean a (continuous) linear functional (possibly 
complex valued)
on the space of continuous functions on $G$ with compact support, which
we identify with a regular Borel measure according to the Riesz-Markov
representation theorem \cite{RS}.

Recall that $G$, when seen as a locally compact space, is $\sigma$-{\em 
compact\/}
(or countable at infinity) if and only if a countable family
\begin{equation} \label{sequence}
     \cA \; =  \; \{ A_n \mid n\in\NN\, \}
\end{equation}
of relatively compact open sets $A_n$ exists with $A_1\neq \es$,
$\overset{\_\!\_}{A_n}\subset A_{n+1}$ for all $n\in\NN$, and
$G=\bigcup_{n\in\NN} A_n$, see \cite[Thm.\ 8.22]{Q}.
In particular, $0<\theta(A_n)<\infty$ for all $n\in\NN$.
This property of $\sigma$-compactness facilitates the calculation
with volume averages which we will need throughout. We assume that
such a sequence $\cA$ has been selected and fixed.

Let $S\subset G$ be a countable set and $w\! : S\to \CC$ a function
such that the corresponding {\em weighted Dirac comb\/}
\begin{equation} \label{comb1}
    \omega \; = \; \sum_{x\in S} w(x)\, \delta_x
\end{equation}
defines a complex regular Borel measure on $G$. Here, $\delta_x$ is
the unit point (or Dirac) measure located at $x$.
The measure $\omega$ need not be bounded, and the unbounded ones are 
those
of particular interest to us here. However, we will assume:
\begin{axiom} \label{ax1}
     {\rm The measure $\omega$ of (\ref{comb1}) is translation bounded.}
\end{axiom}
Recall that a measure $\omega$ is called {\em translation bounded\/}
(or shift-bounded) if, for all compact $K\subset G$,
$\sup_{t\in G} |\omega| (t+K) \le C^{}_K <\infty$ for some constant
$C^{}_K$ which only depends on $K$. Here, $|\omega|$ denotes total
variation measure and $t+K := \{t+x\mid x\in K\}$. If $\omega$ is of
the form (\ref{comb1}), the condition translates into
$\sup_{t\in G} \sum_{x\in S\cap (t+K)} |w(x)| \le C^{}_K <\infty$.
Following \cite{GdeL}, we will denote the space of translation
bounded complex measures on $G$ by $\cM^{\infty}(G)$. For most
of this paper, we consider $\cM^{\infty}(G)$ equipped with the
{\em vague topology}. This means that we identify a translation
bounded measure with the corresponding linear functional on
$\cK (G)$, the space of complex valued continuous functions
with compact support, which is justified by the Riesz-Markov
representation theorem \cite[Thm.~IV.18]{RS}.
Later on, we will also need different topologies on
$\cM^{\infty}(G)$, but we postpone details to Section~\ref{extensions}.

Next, we have to approach the autocorrelation measure attached to 
$\omega$.
With $\omega(g):= \int_G g \, {\rm d}\omega$ for $g\in\cK(G)$,
we can define a partner measure, $\tilde{\omega}$,
via $\tilde{\omega}(g) = \overline{\omega (\tilde{g})} $
where $\tilde{g}(x) := \overline{g(-x)}$.
For the measure $\omega$ of (\ref{comb1}), this means
$\tilde{\omega} = \sum_{x\in S} \overline{w(x)}\, \delta_{-x}$.

Consider the fixed sequence $\cA$ of (\ref{sequence}).
Set $\omega^{}_n = \omega|^{}_{A_n}$ and
$\tilde{\omega}^{}_n = (\omega^{}_n )\tilde{\hphantom{t}}$.
Then, the measure
\begin{equation} \label{auto-approx}
   \gamma^{(n)}_{\omega} \; := \; \frac{\omega^{}_n * 
\tilde{\omega}^{}_n}
   {\theta(A_n)}
\end{equation}
is well defined, since it is the (volume averaged) convolution of two
{\em finite\/} measures. It reads $\gamma^{(n)}_{\omega} =
\sum_{z\in\gD} \eta^{}_n(z) \, \delta_z$ where $\gD = S - S$
(which is still countable) and
\begin{equation}  \label{approx-coeff}
    \eta^{}_n(z) \; = \; \frac{1}{\theta(A_n)}\;
    \sum_{\substack{x,y \in S\cap A_n \\ x-y=z}} w(x)\overline{w(y)}\, .
\end{equation}
For $z\not\in\gD$, we set $\eta^{}_n(z) = 0$. By construction,
$\eta^{}_n$ is then a {\em positive definite function\/} on $G$ (see
\cite[Ch.\ I.3]{BF} for background material). Next,
recall that a measure $\mu$ is called {\em positive definite\/} iff
$\mu(g * \tilde{g})\ge 0$ for all $g\in\cK(G)$. So, the measures
$\gamma^{(n)}_{\omega}$ are positive definite, as follows from
\cite[Prop.\ 4.4]{BF}.
\begin{axiom} \label{ax2}
   {\rm The pointwise limit
       $\eta(z) := \lim_{n\to\infty} \eta^{}_n(z)$
       exists for all $z\in \gD$.}
\end{axiom}
Of course, $\eta(z)=0$ unless $z\in\gD$, so that Axiom \ref{ax2}
guarantees the existence of $\eta(z)$ for all $z\in G$. Also,
$\eta(z)$ is positive definite, and hence satisfies
$|\eta(z)|\le \eta(0)$ for all $z\in G$.

\smallskip \noindent
{\sc Remark}: Let us add a few comments on Axiom \ref{ax2}.
Since $\omega$ is translation
bounded by Axiom \ref{ax1}, the family of finite volume approximations
$\gamma^{(n)}_{\omega}$ of (\ref{auto-approx}) are uniformly
translation bounded \cite[Prop.~2.2]{Hof}. This means that they
always have points of accumulation in the vague topology.
For any such limit point $\gamma$, we can select a suitable subsequence
$\gamma^{(n_i)}_{\omega}$ which converges to $\gamma$.
Since $\gD$ is countable, this can be done in such a way that also
$\eta^{}_{n_i} (z)$ pointwise converges for all $z\in\gD$.
This argument shows that Axiom \ref{ax2} is basically a
convention because we can always achieve it a posteriori by
properly thinning out the (averaging) sequence $\cA$ of 
(\ref{sequence}).

The next Axiom comes in two forms. The first part of the paper requires
only the weaker (first) version of the Axiom. When we come to the
diffraction results, we shall need the stronger version.
\begin{axiom} \label{ax3}
   {\rm The set $\gD$, hence also $\gD^{\rm ess}:=
        \{z\in G\mid \eta(z)\neq 0 \}\subset \gD$, is closed and 
discrete.}
\end{axiom}
This axiom is equivalent to saying that the set $S$ from \eqref{comb1} has 
finite local complexity, compare \cite[Def.~2.2 and Prop.~2.3]{Martin2}.
\addtocounter{axiom}{-1}
\begin{axiom}${\!\!\!}^{\bs{+}}\,$ \label{ax3+}
   {\rm In addition to Axiom 3, the set $\gD^{\rm ess}$ is uniformly 
discrete.}
\end{axiom}
Recall that a set $P\subset G$ is called {\em uniformly discrete\/} if 
an open neighbourhood $V$ of $0\in G$ exists such that $(t+V)\cap P = 
\{t\}$ for all $t\in P$.

After Axiom \ref{ax2} and either Axiom \ref{ax3} or the stronger
Axiom 3$^+$, we see that
the sequence $\big(\gamma^{(n)}_{\omega}\big)_{n\in\NN}$ converges, in 
the
vague topology, to
\begin{equation} \label{autocorr}
   \gamma^{}_{\omega} \; = \;
      \sum_{z\in\gD} \eta(z) \delta_z \; = \;
      \sum_{z\in\gD^{\rm ess}} \eta(z) \delta_z
\end{equation}
which is a pure point measure. In fact, the validity of \eqref{autocorr}
requires Axiom 3. The limit is called the {\em autocorrelation measure\/}, 
or autocorrelation for short, of the Dirac comb $\omega$ of (\ref{comb1}) 
w.r.t.\ the sequence $\cA$ of (\ref{sequence}). It is a translation bounded 
measure, as a consequence of Axiom \ref{ax1}; compare \cite[Prop.~2.2]{Hof}.

Note that different sequences $\cA$ can,
and generally will, lead to different autocorrelations. Nonetheless, the
development of our ideas depends only on the fixed sequence $\cA$.
The effect of changing sequences will be discussed in the Appendix.

Since the positive definite measures are closed in the vague topology,
$\gamma^{}_{\omega}$ is positive definite and $\eta(z)$ is a positive
definite function on all of $G$. Let us assume that the
$\cA$-density of the set $S$ is positive, i.e., $\eta(0)>0$, to
exclude the trivial case that $\gamma^{}_{\omega}=0$.

\smallskip \noindent {\sc Remark}:
In general, $\eta(z)$ will not be a continuous function on $G$, but, by
the Riesz-Segal-von Neumann theorem \cite[p.\ 104]{Gu}, it always admits
a unique decomposition into a continuous positive definite function and
another positive definite function which vanishes $\theta$-almost
everywhere on $G$.
In most cases of interest to us, the continuous part will
be absent, and it is the remaining part that matters.
\smallskip

The support of $\gamma^{}_{\omega}$ (which is the support of $\eta(z)$)
plays a special r\^ole in this article,
and we will in particular need the group generated by it,
\begin{equation} \label{L-def}
     L \; := \; \langle \gD_{}^{\rm ess} \rangle^{}_{\ZZ}
\end{equation}
which is a subgroup of $G$.

Our next goal is the introduction of a
suitable uniformity that allows us to form the (Hausdorff) completion of
$L$. Since $\eta(z)$ is a positive definite function on $G$ (and hence
certainly also on the subgroup $L$), the function
$$  m(z) \; := \; 1 - \frac{\eta(z)}{\eta(0)} $$
is negative definite \cite[Cor.\ 7.7]{BF}. In particular, it satisfies
$m(0) = 0$ and $\tilde{m}=m$, see \cite[Prop.\ 7.5]{BF}. Moreover,
the function $\sqrt{|m(z)|}$ is subadditive, i.e.,
$$ \sqrt{|m(z+z')|} \; \le \; \sqrt{|m(z)|} + \sqrt{|m(z')|} $$
for all $z,z'\in G$, see \cite[Prop.\ 7.15]{BF}. Let us now define
\begin{equation} \label{metric1}
    \varrho(s,t) \; = \; \big| m(s-t) \big|^{1/2}
\end{equation}
which is non-negative, symmetric (due to $\tilde{m}=m$) and satisfies
the triangle inequality (due to the above subadditivity property).
If the weighting function is real and non-negative, one has
$\varrho(s,t)\le 1$ for all $s,t\in G$; in general, $\varrho(s,t)$
is bounded by $\sqrt{2}$.
Also, it is immediate that $\varrho(s+r,t+r) = \varrho(s,t)$ for all
$r\in G$. Hence we have
\begin{fact} \label{fact1}
    The function $\varrho$ of\/ $(\ref{metric1})$ defines a translation
    invariant pseudo-metric, both on $G$ and on $L$.  \qed
\end{fact}

\noindent
{\sc Remark}: If $w(x)\equiv 1$ in the definition (\ref{comb1}) of
$\omega$, one could alternatively work directly with
$\varrho'(s,t):=\eta(0)-\eta(s-t)$. This is then also a pseudo-metric,
called {\em variogram\/} in geostatistics, see \cite[Thm.\ 1]{Matheron}.
However, it is not true in general that negative definite functions
themselves define a pseudo-metric; $g(x)= 1-\exp(-x^2)$ is
a counterexample on $\RR$.

\smallskip
The pseudo-metric $\varrho$ of (\ref{metric1}) defines a uniformity, 
both
on $G$ and on $L$,
and hence a topology, which we call the {\em autocorrelation 
topology\/},
or AC topology for short. This topology is, in general, different from 
the
topology that came with $G$ in the beginning.

Next, we define, for $\ve>0$, a set $P_{\ve}$ of
{\em $\ve$-almost periods\/} of the autocorrelation $\gamma^{}_{\omega}$
through
\begin{equation} \label{almost-periods}
     P_{\ve} \; = \; \{t\in G\mid \varrho(t,0) < \ve \} \, .
\end{equation}
We clearly have $P_{\ve}\subset P_{\ve'}$ for $\ve<\ve'$.
Note that $P_{\ve} = G$ if $\ve > \sqrt{2}$, and that
$\varrho(t,0)=1$ for all $t\in G\setminus\gD_{}^{\rm ess}$.
Consequently, the interesting range of $\ve$ is
$0<\ve\le 1$. Moreover, if $w(x)\ge 0$ on $S$,
$P_1=\gD_{}^{\rm ess}$ and $P_{\ve}=G$
for all $\ve>1$. For $0<\ve\le 1$, we have
$P_{\ve}\subset \gD_{}^{\rm ess}$, because $t\in P_{\ve}$
implies $0 < |\eta(t)| \le \eta(0)$, hence $t\in\gD_{}^{\rm ess}$.
{}From the symmetry
of $\varrho$, we know that $t\in P_{\ve}$ implies
$-t\in P_{\ve}$, while the triangle inequality shows that
$t\in P_{\ve}$ and $t' \in P_{\ve'}$ results in
$t+t' \in P_{\ve+\ve'}$. Thus:
\begin{fact} \label{inclusion}
    The sets $P_{\ve}$ are symmetric, and
    $P_{\ve} + P_{\ve'} \subset
    P_{\ve + \ve'}$. \qed
\end{fact}

\begin{axiom} \label{ax4}
  {\rm For all $\ve>0$, the set $P_{\ve}$ is relatively dense.}
\end{axiom}
Recall that a set $P\subset G$ is called {\em relatively dense\/} if 
there
exists a compact set $K\subset G$ such that, for all $x\in G$,
$(x+K)\cap P \neq \es$. Equivalently, there exists a compact set
$K\subset G$ with $P+K = G$, where $P+K := \{p+k\mid p\in P, \, k\in 
K\}$.
Note that, together with Axiom 3$^+$, we will have
that the $P_{\ve}$ are actually {\em Delone sets\/} if
$0<\ve\le 1$, i.e., they are both uniformly discrete and
relatively dense.

\begin{lemma} \label{krein-bound}
   If $t\in P_{\ve}$, then
   $|\eta(x) - \eta(x+t)| < \eta(0)\, \sqrt{2}\, \ve$,
   for all $x\in G$.
\end{lemma}
{\sc Proof}: If $t\in P_{\ve}$, then
$\big| 1 - \frac{\eta(t)}{\eta(0)}\big|^{1/2} < \ve$ by definition,
which is equivalent to $|\eta(0) - \eta(t)| < \eta(0)\, \ve^2$.
This implies $0 \le \eta(0) - {\rm Re}(\eta(t)) < \eta(0)\, \ve^2$
where the positivity of the middle expression follows from
$|\eta(t)|\le \eta(0)$.

Recall that $\eta$ is a positive definite function on all of $G$, so
for arbitrary $x\in G$, we can invoke Krein's inequality, see
\cite[p.\ 12, Eq.\ (4)]{BF} or \cite[p.\ 103]{Gu}:
$$ \big| \eta(x) - \eta(x+t) \big|^2 \; \le \;
    2\, \eta(0)\, \big( \eta(0) - {\rm Re}(\eta(t)) \big). $$
Together with the previous estimate, this establishes the assertion. 
\qed

\smallskip
The situation which we have outlined has a number of variations.
It is not essential that the original measure $\omega$ be a Dirac
comb. We do require that it be translation bounded and that its
autocorrelation be a point measure supported on a closed discrete set.
Certainly, $\omega$ could have some density $0$ deviation from a
Dirac comb without making any difference to the outcome. In fact, there
is an entire class of measures $\omega$, the so-called homometry class,
all with the same autocorrelation. They cannot be distinguished by
diffraction.

An important type of situation that we wish to point out is that in 
which
$S$ itself is a Delone subset of $G$ with finite local complexity, i.e.,
$S-S$ is closed and discrete. In this case, we can associate with $S$ 
the
Dirac comb $\omega = \delta^{}_S = \sum_{x\in S} \delta_x$,
which obviously satisfies
Axiom~\ref{ax1}. If it also satisfies Axiom~\ref{ax2}, we say that it
has a well-defined autocorrelation (relative to the sequence $\cA$).
Axiom~\ref{ax3} is then already
implicit in the finite local complexity, which has the geometric meaning
that for any compact set $K \subset G$ there are only finitely many 
types
of point sets $S\cap (a+K)$, up to translation, as $a$ runs over $G$.
Furthermore, in this case, it is possible to replace $\varrho$ by the
topologically equivalent but simpler pseudo-metric $\varrho'$
which was mentioned in the Remark after Fact~\ref{fact1}.
We also say, by slight abuse of language, that $S$ is $\ve$-almost
periodic if Axiom~\ref{ax4} holds. The results of Section 3 below
apply to point sets $S$ which satisfy these conditions.

Let us close this part with a basic result on the diffraction measure
attached to an arbitrary $\omega\in\cM^{\infty}(G)$ for a specified
sequence $\cA$ of (\ref{sequence}). If the autocorrelation measure
$\gamma^{}_{\omega}$ exists as the limit of the $\gamma^{(n)}_{\omega}$
defined in (\ref{auto-approx}), then 
$\gamma^{}_{\omega}\in\cM^{\infty}(G)$,
see \cite[Prop.\ 2.2]{Hof} and our above Remark after Axiom~\ref{ax2}.
It is also positive definite because the positive definite measures
are closed in the vague topology \cite[p.\ 18]{BF}. As such, the
measure $\gamma^{}_{\omega}$ is transformable, i.e., there is a
uniquely determined positive measure $\hat{\gamma}^{}_{\omega}$
on $\hat{G}$ such that
$$ \gamma^{}_{\omega} \big( g*\tilde{g} \big)  \; = \;
    \hat{\gamma}^{}_{\omega} \big( |\tilde{\hat{g}}|^2 \big)$$
for all $g\in\cK(G)$, where $\hat{g}$ denotes ordinary Fourier
transform on $L^1(G)$, compare
\cite[Thm.\ 4.7 and Cor.\ 4.8]{BF}. What is more,
$\hat{\gamma}^{}_{\omega}$ is also translation bounded
\cite[Prop.\ 4.9]{BF}. To summarize:
\begin{fact} \label{general-fourier}
   Let $G$ be a $\sigma$-compact LCA group, and 
$\omega\in\cM^{\infty}(G)$.
   Let an averaging sequence $\cA$ be fixed.
   If the corresponding autocorrelation
   $\gamma^{}_{\omega}$ exists, it is a positive definite, translation
   bounded measure on $G$. Moreover, its Fourier transform
   $\hat{\gamma}^{}_{\omega}$ exists and is a positive,
   translation bounded measure on $\hat{G}$, the dual group of $G$.
   \qed
\end{fact}

\section{Construction of a cut and project scheme}\label{sec2}

The goal of this section is to establish how $L$ is related, in a 
natural
way, to a lattice in an appropriate extension of the LCA group $G$, and 
to use
this connection to construct a cut and project scheme that will 
ultimately
help us to prove the pure point diffractivity. Our assumptions in this
Section are Axioms 1 -- 4, but not Axiom 3$^+$.

An essential ingredient of our approach is the precompactness of certain
sets. Recall that a subset $P$ of the LCA group $G$ is {\em 
precompact\/}
(or totally bounded)
if, for any open neighbourhood $U$ of $0\in G$, there is a {\em 
finite\/} set $F$
so that $P\subset U+F$. Note that $P$ is precompact if and only if the
closure of its image in the Hausdorff completion is compact, compare
\cite[Ch.\ 13 A]{Q} and \cite[Ch.\ II.4.2]{Bou}.

\begin{lemma} \label{precomp1}
    For all\/ $0<\ve < 1$, the set $P_{\ve}$ is precompact in the AC 
topology.
\end{lemma}
{\sc Proof}: In view of the above remark, it suffices to show that,
for any $\ve'>0$, $P_{\ve}$ can be covered by finitely
many translates of the set $P_{\ve'}$ because the latter form a 
fundamental
system of neighbourhoods.
Let $0<\ve <1$ be fixed, and choose any $\ve'$ with
$0 < \ve' < 1-\ve$. Since  $P_{\ve'}$ is
relatively dense, there exists a compact set $K\subset G$ such that
$P_{\ve'} + K = G$.

Let $t\in P_{\ve}$ and write $t=s+k$ with $s\in P_{\ve'}$
and $k\in K$. Then we get
$$ t - s \; \in \; P_{\ve} - P_{\ve'} \; = \;
    P_{\ve} + P_{\ve'} \; \subset \;
    P_{\ve+\ve'} $$
from Fact \ref{inclusion}.
We conclude that $k=t-s \in F := ( K\cap P_{\ve+\ve'}) \subset L$ which 
is
a {\em finite\/} set because $\ve + \ve' < 1$ by construction, so we
know that $P_{\ve+\ve'}\subset \gD_{}^{\rm ess}$ which is
discrete and closed.
So, $t\in P_{\ve'} + F$, and
hence $P_{\ve}\subset P_{\ve'}+F$. This being true
for all $0 < \ve' < 1-\ve$, it is clearly true for
all $\ve'>0$ and our assertion follows.   \qed

\smallskip
Note that the proof does {\em not\/} extend to $\ve=1$.
In fact, we presently see no general reason for
$P^{}_1 = \gD_{}^{\rm ess}$ to be precompact. However,
it is precompact if $S$ (the support of $\omega$) is a
model set, or the set of visible lattice points.

\smallskip
Let us now define $H$ to be the (Hausdorff) completion of $L$ of 
(\ref{L-def})
in the AC topology. We then know that a uniformly continuous
homomorphism $\varphi\! : L\to H$ exists with the following properties:
\begin{itemize}
\item [{\bf C1}] The image $\varphi(L)$ is dense in $H$.
\item [{\bf C2}] The mapping $\varphi$ is an open mapping from $L$ onto
       $\varphi(L)$, the latter with the induced topology of the 
completion $H$.
\item [{\bf C3}] We have $\ker(\varphi) =
         \mbox{closure of } \{0\}\mbox{ in } L$.
\end{itemize}
In fact, these three properties characterize $(H,\varphi)$ as {\em 
the\/}
completion of $L$, in the sense that it is unique up to isomorphism.
Note that $\ker(\varphi) = \bigcap_{\ve>0} P_{\ve}$ and that
the condition for $\varphi$ to be an embedding is thus
$$ \bigcap_{\ve>0} P_{\ve} \; = \; \{0\} \, . $$

Since this characterization of the completion is not entirely standard, 
we
add a few comments. Since the Hausdorff completion is unique up to
isomorphism, the three properties follow rather directly from the 
universal
characterization in \cite[Thm.~II.3.3]{Bou} (recall that a bijection 
$\pi$ between
two uniform spaces is an isomorphism if and only if both $\pi$ and its
inverse are uniformly continuous \cite[Prop.~I.2.2 (b)]{Bou}).
Conversely, given a homomorphism $\varphi$
with the three properties, one compares with the Hausdorff completion, 
which
now brings in another mapping, $\psi\!: L\to \hat{L}$ say. Then, the
existence of a uniformly continuous mapping $\pi\!: \hat{L}\to H$ with
$\pi\circ\psi=\varphi$ is assured. Comparing kernels, one concludes that
$\pi$ is one-to-one on $\psi(L)$. So, with the open mapping property of
$\varphi$, we may conclude that the dense subsets $\varphi(L)$ and
$\psi(L)$ are isomorphic (via $\pi$), and this then extends to an
isomorphism between $H$ and $\hat{L}$ by the Corollary of
\cite[Thm.~II.3.2]{Bou}.

\begin{prop} \label{HisLCAG}
    The completion $H$ is a locally compact Abelian group.
\end{prop}
{\sc Proof}: That $H$ is a topological group and Abelian is
an immediate consequence of the completion procedure. For the
local compactness, it is sufficient to show that $0\in H$ has a
compact neighbourhood. Fix some $0<\ve<1$. The set $P_{\ve}$
is precompact (Lemma \ref{precomp1}), hence
$\overline{\varphi(P_{\ve})}\subset H$ is compact, compare
\cite[Thm.\ 13.2]{Q}. It is, in fact, a neighbourhood of $0\in H$.

To see this, observe that $\varphi(P_{\ve}) = \varphi(L)\cap V$
for some open set $V\subset H$ because of property C2.
Since $0\in P_{\ve}$, $\varphi(P_{\ve})$ contains $0\in H$,
so we also have $0\in V$, and $V$ is an open neighbourhood of $0$.
If $v\in \overline{V}$, for every open neighbourhood $N$ of $v$,
$N\cap V\neq\es$, hence also $N\cap V\cap\varphi(L)=
N\cap \varphi(P_{\ve})\neq \es$.
Consequently, $v\in\overline{\varphi(P_{\ve})}$ and
$\overline{V}=\overline{\varphi(P_{\ve})}$.  \qed

\smallskip
Since $\varphi(L)\subset H$, we now also have a mapping from $L$ to
$G\times H$, sending $t\in L$ to the point $(t,\varphi(t))$. Let
$\tilde{L}$ be the image of $L$ under this mapping.
\begin{lemma} \label{lattice1}
    The set $\tilde{L}$ is uniformly discrete in $G\times H$.
\end{lemma}
{\sc Proof}: Fix some $\ve$ with $0<\ve<1$ and choose
an open set $V\subset H$ such that
$V\cap\varphi(L)=\varphi(P_{\ve})$, which is always
possible due to property C2.
Let $U$ be a compact neighbourhood of $0\in G$ and assume
that $(t,\varphi(t))\in U\times V$ for some $t\in L$.
Then $\varphi(t)\in V\cap\varphi(L)=\varphi(P_{\ve})$
which implies $t\in P_{\ve}$ (since $P_{\ve}
\supset \ker(\varphi) \supset \{0\}$ and
$P_{\ve} + \ker(\varphi) = P_{\ve}$),
hence also $t\in U\cap P_{\ve}=F$ where $F$ is a finite set.
Consequently, $\tilde{L}\cap (U\times V)$ is finite, too.

Since $U\times V$ is a neighbourhood of $(0,0)\in G\times H$,
$(0,0)$ is isolated in $\tilde{L}$, and by translation this
is true of all points of $\tilde{L}$. The Hausdorff property of $H$
then implies that $\tilde{L}$ is discrete, with a uniform separating
neighbourhood, hence also uniformly discrete. \qed

\smallskip
{}For any subset $X\subset H$, we define
\begin{equation} \label{model-def}
    \oplam(X) \; := \; \{ x\in L\mid \varphi(x)\in X \} \, .
\end{equation}
This type of construction of subsets in $G$ plays an important
role in the sequel.

\begin{lemma} \label{model1}
    {}For all\/ $0<\ve<1$, the compact set\/ 
$W:=\overline{\varphi(P_{\ve})}$
    has non-empty interior, i.e., $\overset{\circ}{W}\neq\es$.
\end{lemma}
{\sc Proof}: It follows from the proof of Proposition~\ref{HisLCAG}
that $\overline{\varphi(P_{\ve})}$ is a (closed) neighbourhood of
$0\in H$, and hence has non-empty interior.  \qed

\smallskip
Let us briefly detour to show an alternative argument for the
validity of Lemma~\ref{model1} on the basis of Baire's theorem.

 From Lemma \ref{precomp1}, $P_{\ve}$ is precompact
(in $L$), so $W$, and also $\partial W$, is compact in the Hausdorff
completion. Then $\oplam(W)$ clearly
contains $P_{\ve}$ and is thus relatively dense in $G$. The strategy
will now be to show that $\oplam(\partial W)$ is {\em not\/} relatively 
dense
in $G$, hence $W\neq\partial W$, and $\overset{\circ}{W}\neq\es$.

We employ Baire's theorem \cite[Thm.\ 13.29 (b)]{Q} which applies to
locally compact spaces, so $H$ is of second category as a consequence
of Proposition~\ref{HisLCAG}. Recall that $\gD=S-S$ is countable, so 
also
$L=\langle \gD_{}^{\rm ess} \rangle^{}_{\ZZ}$ is still countable, and
$\varphi(L)$ is then a meagre subset of $H$.
Since $\partial W$ is nowhere dense, we know
that a $c\in H$ exists such that $(c+\partial W)\cap\varphi(L) = \es$.
Let $K\subset G$ be any non-empty compact set. We then clearly have
$\big(K\times (c+\partial W)\big) \cap \tilde{L} = \es$.

If $V_1$ is a compact neighbourhood of $0$, then
$K\times (c+\partial W + V_1)$ is compact and
$\big(K\times(c+\partial W + V_1)\big)\cap\tilde{L}$ is finite.
We can then find a neighbourhood $V$ of $0\in H$
such that $\big(K\times (c+\partial W+V)\big)\cap \tilde{L} = \es$.

{}Finally, since $\varphi(L)$ is dense in $H$, there exists $t\in L$ 
with
$-\varphi(t) \in c+V$, so that
$\big(K\times (-\varphi(t) + \partial W)\big)\cap\tilde{L} = \es$.
This implies $\big((t+K)\times \partial W\big)\cap\tilde{L} = \es$
and hence $(t+K)\cap\oplam(\partial W) = \es$, where $t$ depends on $K$.
Since the compact set $K\subset G$ was arbitrary, this means
$\oplam(\partial W)$  is not relatively dense.

\begin{lemma} \label{lattice2}
    The set $\tilde{L}$ is relatively dense in $G\times H$.
\end{lemma}
{\sc Proof}: Fix some $0<\ve<1$ and set
$W=\overline{\varphi(P_{\ve})}$ as in Lemma \ref{model1}.
Since $P_{\ve}\subset L$
is relatively dense in $G$ by Axiom \ref{ax4}, there is a compact set
$K\subset G$ such that $G=P_{\ve} + K$. Consequently,
$G\times \{0\} \subset \tilde{L} + \big(K\times (-W)\big)$, because
any $x\in G$ now admits the decomposition
$$ (x,0) \; = \; (t+k,0) \; = \; (t,\varphi(t)) + (k,-\varphi(t)) $$
for some $t\in P_{\ve}$ and $k\in K$, where
$(t,\varphi(t))\in\tilde{L}$.

{}Furthermore,
$G\times H = \tilde{L} + (G\times W)$ since $\varphi(L)$ is dense in
$H$ by property C1 and $\overset{\circ}{W}\neq\es$ by Lemma 
\ref{model1}.
With $G\times W = (G\times\{0\}) + (\{0\}\times W)$, we get
\begin{eqnarray*}
    G\times H & = & \tilde{L} + (G\times W) \;\, \subset \;\,
    \tilde{L} + \tilde{L} + \big(K\times (-W)\big) + (\{0\}\times W) \\
    & = & \tilde{L} + \big(K\times (W-W)\big) \;\, \subset \;\,
    G\times H .
\end{eqnarray*}
Consequently, $G\times H = \tilde{L} + \big(K\times (W-W)\big)$,
and since $K\times (W-W)$ is compact, $\tilde{L}$ is
relatively dense.  \qed

\smallskip
So, $\tilde{L}$ is a discrete subgroup of $H$, hence closed, and it
is also relatively dense, hence co-compact. Thus:
\begin{coro}
    The factor group\/ $\TT:=(G\times H) / \tilde{L}$ is compact.  \qed
\end{coro}

In what follows, we give $G$ the structure of a topological group in
a new way that mixes its standard topology with the AC topology on $L$.
The main result will be that the completion of $G$ with respect to
this new topology is precisely the group $\TT$.

Consider $G\times L$ with the product topology, where $G$ is given the
standard topology it came with and $L$ is given the AC topology that was
introduced in Section \ref{sec1}. Let
$$  \alpha : L \longrightarrow G\times L $$
be the diagonal map, i.e., $\alpha(t) = (t,t)$. For any open 
neighbourhood
$U$ of $0\in G$ and any $0<\ve \le 1$, we then have
$$ \alpha(L) \cap (U\times P_{\ve}) \; = \;
    \{ (t,t) \mid t \in U\cap P_{\ve} \} $$
which is finite. Thus, $\alpha(0)$ is isolated in $G\times L$ and
$\alpha(L)$ is a closed discrete subgroup of $G\times L$.

Consequently, $(G\times L)/\alpha(L)$ is a topological group with the
natural quotient topology. The natural map
\[
    \psi \! : \; G\times L \longrightarrow (G\times L)/\alpha(L) 
\]
is continuous and also open. Furthermore, $G\simeq (G\times L)/\alpha(L)$
through the mapping $x \mapsto (x,0)$ mod $\alpha(L)$, and in this way 
we obtain a new `mixed' topology on $G$. With $\varphi \!: L\to H$ as
before, we have the following diagram:
\[
\begin{CD}
    0 @>>> L @>\alpha>> G\times L    @>\psi>> (G\times L)/\alpha(L) @>>> 0 \\
    @. @V{\rm id}VV @V{{\rm id}\times\varphi}VV @VV\overline{\varphi}V @. \\
    0 @>>> L @>{\rm id}\times\varphi>> G\times H @>\pi>>
           \TT=(G\times H)/\tilde{L} @>>> 0
\end{CD}
\]
Here, $\pi$ is the natural map from $G\times H$ to $\TT$, so that the 
two
horizontal sequences are exact, and $\overline{\varphi}$ is the unique
(continuous) homomorphism which makes the central block of the diagram
commutative.

\begin{prop} \label{completion1}
    The mapping\/ $\overline{\varphi}\! : (G\times L)/\alpha(L)
    \longrightarrow \TT$ defines a completion of $(G\times L)/\alpha(L)$.
    In particular, $\TT$ may be regarded as the Hausdorff completion of
    $G$ with respect to the mixed topology.
\end{prop}
{\sc Proof}: Since $\TT$ is a compact group, it is complete
\cite[Cor.\ II.3.1]{Bou}. It suffices to establish the three
properties (C1)--(C3) for $(G\times L)/\alpha(L)$, $\TT$, and
$\overline{\varphi}$.

(I) ${\rm im}(\overline{\varphi})$ is dense in $\TT$:  \newline
The closure of the image of $G\times L$ in $\TT$ has a closed
preimage in $G\times H$ which contains $G\times \varphi(L)$ and
hence also $G\times H$.

(II) $\ker(\overline{\varphi})$ is the closure of $\{0\}$ in
     $(G\times L)/\alpha(L)$:  \newline
Certainly, $\ker(\overline{\varphi})$ is closed. Let
$\psi(x,t)\in\ker(\overline{\varphi})$. This implies
$(x,\varphi(t))\in\tilde{L}$, so $x\in L$ and
$\varphi(x) = \varphi(t)$. Then, $(t-x) \in \ker(\varphi)$ and
$\psi(x,t) = \psi(0,t-x)$ is contained in
$\psi\big(\{0\}\times\ker(\varphi)\big)$.
Since $\ker(\varphi)$ is the closure of $\{0\}$ in
$L$, each neighbourhood of $(t-x)$ in $L$ contains $0$.
So, each neighbourhood of $\psi(x,t)$ must contain $\psi(0,0)$.

(III) $\overline{\varphi}$ is an open mapping onto
    ${\rm im}(\overline{\varphi})$,
    the latter equipped with the topology induced from $\TT$:  \newline
Let $U\subset G$, $V\subset L$ be non-empty open sets. Then
$\varphi(V) = W\cap\varphi(L)$ for some open set $W\subset H$
and
$$\begin{CD}
   U\times V  @>\psi>> \psi(U\times V) \\
   @V{{\rm id}\times\varphi}VV @VV\overline{\varphi}V \\
   \begin{array}{c} U\times\varphi(V) = \\
    (U\times W)\cap(G\times\varphi(L)) \end{array}
    @>\pi>> \begin{array}{c}
    \pi\big( (U\times W)\cap(G\times\varphi(L)) \big) = \\
    \pi(U\times W) \cap \pi(G\times\varphi(L))
    \end{array}
\end{CD}$$
The only part of this diagram which needs explanation is the equation
in the bottom right hand corner. With the obvious notation, one gets
\begin{eqnarray*}
   \lefteqn{\pi(u,w) \; = \; \pi(x,\varphi(t))
     \qquad \text{(for some $x\in G$ and $t\in L$)} } \\
   & \Longrightarrow &
     (u,w) \; = \; (x,\varphi(t)) + (s, \varphi(s))
     \qquad \text{(for some $s\in L$)} \\
   & \Longrightarrow &
     (u,w) \; = \;  (x+s,\varphi(t+s))
     \; \in \; (U\times W)\cap \big(G\times \varphi(L)\big).
\end{eqnarray*}
This completes the proof.  \qed

\smallskip
\noindent
{\sc Remark}: In plainer language, the new mixed topology on $G$
can be described by saying that $x,y\in G$ are close if, for some
small $v\in G$ and some $t\in P_{\ve}$ with small $\ve>0$, we have
the relation $y=x+v+t$.

\section{Diffraction}

In this Section, we will be assuming Axioms \ref{ax1}, \ref{ax2},
3$^+$, and \ref{ax4}.
In particular, the constructions and results of Section \ref{sec2}
are in force.

The plan is now to use the relation between $G$ and $\TT$ to derive the
nature of the Fourier transform of the autocorrelation 
$\gamma^{}_{\omega}$.
This will be achieved by first regularizing it through the convolution 
with
suitable continuous functions, then studying the properties of their 
Fourier
transforms, and later transfering them to $\hat{\gamma}^{}_{\omega}$ via
appropriate limits.

So, let $c\in \cK(G)$ be a real valued and non-negative function.
Then, also $c * \tilde{c} \in \cK(G)$, i.e., it has compact support.
Consequently,
$g^{}_c = (c*\tilde{c})*\gamma^{}_{\omega}$ is a well-defined
bounded continuous function, see \cite[Prop.\ 1.12]{BF}, because
$\gamma^{}_{\omega}$ is translation bounded.
It is also a positive definite function by construction,
hence transformable, i.e., its Fourier transform exists and is a
finite positive measure on the dual group, $\hat{G}$, by
Bochner's theorem \cite[Thm.\ 3.12]{BF}.
Using the convolution theorem \cite[Prop.\ 4.10]{BF}, we see that
$$ \hat{g}^{}_c \; = \; |\hat{c}|^2 \, \hat{\gamma}^{}_{\omega}$$
where $\hat{c}$ is a continuous function on $\hat{G}$.

Let us find out more about $\hat{g}^{}_c$. Observe first that $g^{}_c$
is a positive definite function which is real and certainly continuous
at $0\in G$. By \cite[Prop.\ 3.10]{BF}, it is then automatically a
uniformly continuous function on $G$ in the original topology.
But in view of our assumptions, we can show more:
\begin{lemma} \label{cont1}
      The function  $g^{}_c$ is uniformly continuous in the mixed 
topology.
\end{lemma}
{\sc Proof}: It is sufficient to show that, for all $0 <\ve \le 1$,
$t\in P_{\ve}$ implies the estimate $|g^{}_c(x) - g^{}_c(x+t)| \le C 
\ve$
for all $x\in G$ and for some $C=C(c)$. A direct calculation shows
$$ g^{}_c (x) \; = \; \sum_{z\in\gD_{}^{\rm ess}} \eta(z)\,
    \big(c * \tilde{c}\big) (x-z) \; = \;
    \sum_{z\in L} \eta(z)\, \big(c * \tilde{c}\big) (x-z) \, . $$
Consequently, for $t\in P_{\ve}$, where $P_{\ve}\subset L$, one obtains
\begin{eqnarray*}
    \big| g^{}_c (x) - g^{}_c (x+t) \big| & = &
    \Big| \sum_{z\in L} \big( \eta(z) - \eta(z+t) \big)\,
    \big(c * \tilde{c}\big) (x-z) \Big|  \\
    & \le & \| c * \tilde{c} \|^{}_{\infty}
      \sum_{z\in F_x} \big| \eta(z) - \eta(z+t) \big| \\
    & \le &  \| c * \tilde{c} \|^{}_{\infty} \, \sqrt{2} \,
      \eta(0)\, {\rm card}(F_x)\, \ve
\end{eqnarray*}
where Lemma \ref{krein-bound} was used in the last step and $F_x$ is 
given by
$$F_x = \big( \gD_{}^{\rm ess} \cap (x - {\rm supp}(c*\tilde{c}))\big)
   \cup \big( (-t + \gD_{}^{\rm ess}) \cap (x - {\rm 
supp}(c*\tilde{c}))\big).$$
This is a finite set because $\gD_{}^{\rm ess} $ is closed and discrete 
by
Axiom \ref{ax3} and $(c*\tilde{c})$ has compact support.
In fact, by Axiom 3$^+$,
the cardinality of $F_x$ is uniformly bounded in $x$ by a constant
that is independent of $t$, so our above estimate is uniform in $x$ and 
the assertion follows.  \qed

\smallskip
The purpose of this exercise is that we can now extend $g^{}_c$ to a 
uniquely
defined uniformly continuous function on $\TT$, the completion of $G$ 
with the
mixed topology. We will call this new function $g^{\TT}_c$.

As the uniformly continuous extension of $g^{}_c$, also $g^{\TT}_c$ is a
continuous positive definite function, this time on all of $\TT$. By
Bochner's theorem, see \cite[Thm.\ 3.12]{BF} or \cite[Thm.\ 
1.4.3]{Rudin},
there is a uniquely defined
finite positive measure $\mu^{}_c$ on the dual group, $\hat{\TT}$, so 
that
$$ g^{\TT}_c(x) \; = \;
       \int_{\hat{\TT}} \langle k,x\rangle \, {\rm d}\mu^{}_c (k) $$
for all $x\in G$. Here, $\langle k,\cdot\rangle$ denotes the continuous
character on $\TT$ attached to $k\in\hat{\TT}$.
Furthermore, since $\TT$ is compact, $\hat{\TT}$ is discrete, and we 
thus
get a representation of $g^{\TT}_c$ as an absolutely convergent Fourier
series, i.e.,
\begin{equation} \label{fourier1}
   g^{\TT}_c (x) \; = \; \sum_{k\in \hat{\TT}} a^{}_c (k)\, \langle 
k,x\rangle
\end{equation}
where equality holds for all $x\in\TT$, compare \cite[Thm.\ 2.8.4 
(ii)]{Edwards}.
Here, the $a^{}_c (k)$ are the non-negative numbers
$$ a^{}_c (k) \; = \; \int_{\TT} \overline{\langle k,x\rangle} \,
    g^{\TT}_c(x) \, {\rm d}\theta^{}_{\TT} (x) \; = \;
    \mu^{}_c (\{k\}) \; \ge \; 0  $$
with $\sum_{k\in\hat{\TT}}\, a^{}_c (k) = g^{\TT}_c (0) < \infty$.

We have already proved that
the canonical homomorphism $\beta\! : G \rightarrow \TT =
(G\times H)/\tilde{L}$ is given by
$$ x \; \mapsto \; (x,0) \mbox{ mod }\tilde{L}\, . $$
This leads to a mapping  $C(\TT,\CC) \rightarrow C(G,\CC)$ where
$h^{\TT}\mapsto h^{\TT}\circ\beta$ and, in particular,
to an injective mapping
$\hat{\beta}\! : \hat{\TT}\simeq \tilde{L}^{\circ} \rightarrow \hat{G}$.
Here, $\tilde{L}^{\circ}$ is the dual of the lattice $\tilde{L}$, i.e.,
the space of all continuous homomorphisms of $G\times H \rightarrow \CC$
that are trivial on $\tilde{L}$. Since $\beta(G)$ is dense in $\TT$,
$\hat{\beta}$ is one-to-one.

{}For $k\in\hat{\TT}$, we write $\ell^{}_k$ for the corresponding 
character
on $G$ attached to $\hat{\beta}(k)$, i.e.,
$ \ell^{}_k (x) = \langle k, \beta(x) \rangle$.
Thus, with (\ref{fourier1}), we now obtain for the function $g^{}_c$ the
expansion
$$ g^{}_c (x) \; = \; \sum_{k\in\hat{\TT}} a^{}_c(k) \,
     \ell^{}_k (x) $$
which is convergent for all $x\in G$. Consequently, for all
$h\in \cK(\hat{G})$,
\begin{eqnarray*}
   \big( \hat{g}^{}_c , h \big)  & = &  \big( g^{}_c , \hat{h} \big)
   \;\, = \;\, \int_G \hat{h}(x) \sum_{k\in\hat{\TT}} a^{}_c(k)\,
               \ell^{}_k (x) \, {\rm d}\theta^{}_G (x) \\
   & = & \sum_{k\in\hat{\TT}} a^{}_c(k) \int_G \hat{h}(x) \,
         \ell^{}_k (x) \, {\rm d}\theta^{}_G (x)
   \;\, = \;\, \sum_{k\in\hat{\TT}} a^{}_c(k) \, h(\hat{\beta}(k)) 
\\[1mm]
   & = & \Big( \sum_{k\in\hat{\TT}} a^{}_c(k)\,
         \delta^{}_{\hat{\beta}(k)} , h  \Big)
\end{eqnarray*}
where the interchange of the sum and the integral is justified by the
summability of $\sum_{k\in\hat{\TT}} a^{}_c(k)$. This shows that
$\hat{g}^{}_c$ is a positive pure point measure,
\begin{equation} \label{fourier2}
    \hat{g}^{}_c \; = \; \sum_{\ell \in L_{}^{\circ}}
    a^{}_c \big(\hat{\beta}^{-1} (\ell)\big) \, \delta^{}_{\ell} \, ,
\end{equation}
with $L_{}^{\circ} = \hat{\beta}(\tilde{L}^{\circ})$.

The next step consists in choosing a net of functions
$(c^{}_{U})^{}_{U\in\,{\mathcal U}}$ that constitutes an
{\em approximate unit\/}, compare \cite[Def.\ 1.6]{BF}, where
$\mathcal U$ is a basis for the neighbourhood filter of $0\in G$.

So, let $c^{}_U\in \cK(G)$ be a non-negative function, with
${\rm supp}(c^{}_U) \subset U$ and normalization
according to $\int_{G} c^{}_U(x) \,{\rm d}\theta^{}_G = 1$,
and let the net of probability measures $c^{}_U \, \theta^{}_G$
converge to $\delta_0$ in the vague topology. Such a net exists
for all LCA groups. Then, also
$(c^{}_U * \tilde{c}^{}_U)\, \theta^{}_G \to \delta_0$ vaguely,
and the `forward version' of
Levy's continuity theorem \cite[Thm.\ 3.13]{BF} gives us compact 
convergence
of $|\hat{c}^{}_U|^2$ towards $\hat{\delta}_0 \equiv 1$ on $\hat{G}$.
For simplicity of notation, we now write
$$ \gamma^{}_U \; := \; g^{}_{c^{}_U} \; = \;
    c^{}_U * \tilde{c}^{}_U * \gamma^{}_{\omega} \, .$$
By what we have just seen, $\hat{\gamma}^{}_U$ is a positive pure
point measure on $\hat{G}$, supported on $L^{\circ}$. Also,
$\hat{\gamma}^{}_U = |\hat{c}^{}_U|^2 \, \hat{\gamma}^{}_{\omega}$
by the convolution theorem.
\begin{prop} \label{reg-conv}
    Let\/ $(c^{}_U)^{}_{U\in\,{\mathcal U}}$ be a normalized regularization net
    of\/ $\delta_0$ as described above, and let $\omega$ be the Dirac
    comb of\/ $(\ref{comb1})$. Then, the measures
    $\hat{\gamma}^{}_U$ converge to $\hat{\gamma}^{}_{\omega}$,
    both in the vague and in the
    $K\!$-local norm topology,\footnote{For any
    compact $K\subset \hat{G}$, one has
    $\|(\hat{\gamma}^{}_U - \hat{\gamma}^{}_{\omega})|_K\|
    \to 0$, where $\|.\|$ is the usual variation norm for {\em finite\/}
    measures.}
     for any compact set $K$.
\end{prop}
{\sc Proof}: The vague convergence is a direct consequence of the
continuity of the Fourier transform. To show the
second assertion, we have to demonstrate that the restriction of
$\hat{\gamma}^{}_U$ to any compact set $K\subset \hat{G}$ norm converges
to the restriction of $\hat{\gamma}^{}_{\omega}$ to $K$.

Let $K\subset\hat{G}$ be compact and fix $\ve > 0$. From the above,
we know that
$\sup_{x\in K}\, \big| |\hat{c}^{}_U (x)|^2 - 1 \big| < \ve$,
for all sufficiently small neighbourhoods $U$ of $0\in G$.
Since $\hat{\gamma}^{}_{\omega}$ is a positive measure,
we then also have
$$ \bigl\| \big( \hat{\gamma}^{}_U - \hat{\gamma}^{}_{\omega}
    \big)\big|^{}_K \big. \bigr\| \; = \; \int_K \big| |\hat{c}^{}_U|^2 - 
1 \big|
    \, {\rm d}\hat{\gamma}^{}_{\omega}
    \; \le \; \int_K \, \sup_{x\in K}\, \big| |\hat{c}^{}_U(x)|^2 - 1 
\big|
    \, {\rm d}\hat{\gamma}^{}_{\omega}
    \; \le \; \ve \, \hat{\gamma}^{}_{\omega} (K)
$$
from which the assertion follows.  \qed

\smallskip
Proposition \ref{reg-conv} is an important cornerstone because both
the absolutely continuous and the pure point measures are norm closed
subsets in the cone of positive measures, compare \cite[Ch.\ IV.5]{RS}.
This rests upon the equation
$|\nu^{}_1 + \nu^{}_2| = |\nu^{}_1| + |\nu^{}_2|$ which is
valid if $\nu^{}_1 \perp \nu^{}_2$, i.e., for measures
which are mutually singular. In particular, this applies if
$\nu^{}_1$ is pure point and $\nu^{}_2$ continuous.
If a sequence of pure point measures $\mu_n$ norm converges
to $\mu = \mu_{\rm pp} + \mu_{\rm cont}$, one
can then show that, for any compact set $K$ and any $\ve>0$,
$|\mu_{\rm cont}|(K) < \ve$, hence $\mu_{\rm cont}=0$.
By construction, the $\hat{\gamma}^{}_n$ are all pure point measures
and form a norm-converging sequence. The limit
$\hat{\gamma}^{}_{\omega}$ is thus also a positive pure point measure,
and translation bounded by Fact~\ref{general-fourier}.
To summarize:

\begin{theorem}  \label{main-thm}
    Let $G$ be a $\sigma$-compact LCA group.
    If the Dirac comb\/ $\omega$ of\/ $(\ref{comb1})$, seen as a
    complex measure on $G$, and its autocorrelation
    $\gamma^{}_{\omega}$ relative to the averaging sequence $\cA$ of\/
    $(\ref{sequence})$ satisfy Axioms\/ {\rm 1}, {\rm 2}, {\rm 3$^+$},
    and {\rm 4}, then
    the corresponding diffraction measure $\hat{\gamma}^{}_{\omega}$
    on the dual group $\hat{G}$ is a translation bounded,
    positive pure point measure. \qed
\end{theorem}

\section{Diffraction in model sets}\label{sec4}

A {\em cut and project scheme} is a triple $(G,J, \tilde M)$
consisting of a pair of LCA groups $G,J$
and a lattice $\tilde M \subset G \times J$ for which the
canonical projections $\pi^{}_G : G \times J \longrightarrow G$ and
$\pi^{}_J : G \times J \longrightarrow J$ satisfy
\begin{enumerate}
\item $\pi^{}_G|_{\tilde M}$ is one-to-one and
\item $\pi^{}_J(\tilde M)$ is dense in $J$.
\end{enumerate}

We write $M:= \pi^{}_G(\tilde M)$, so $M$ is a subgroup
of $G$ and note that the mapping
\begin{equation} \label{star}
    (\;)^* \; := \; \pi^{}_J \circ (\pi^{}_G|_{\tilde M})^{-1} :
           \; M \longrightarrow J
\end{equation}
has dense image in $J$. Note that, for $x\in\tilde{M}$,
$(\pi^{}_G|_{\tilde M})^{-1}(x)
= (\pi^{}_G)^{-1} (x) \cap \tilde{M}$ which is a single point by
our assumption on $\pi^{}_G$. We will denote by $\theta_G$ and 
$\theta_J$
a fixed pair of Haar measures on $G$ and $J$.

Let us summarize our findings of Section~\ref{sec2} for later use.
\begin{prop}
    Let $G,J$ be LCA groups, and assume that $G$
    is also $\sigma$-compact. Then $(G,H,\tilde{L})$
    with $\tilde{L}$ as constructed in Section~\ref{sec2} constitutes a
    cut and project scheme, where $H$ is the Hausdorff completion of
    $L=\langle \gD_{}^{\rm ess} \rangle_{\ZZ}$ in the AC topology.
    Moreover, $\pi^{}_G (\tilde{L}) = L$  and $\varphi$ is the
    $*$-mapping  of\/ $(\ref{star})$, hence $\pi^{}_H (\tilde{L})= 
\varphi(L)$.
\qed
\end{prop}

A set $\gL \subset G$ is a {\em regular model set\/} for the cut
and project scheme $(G,J,\tilde{M})$ if there is a relatively
compact set $W \subset J$ with non-empty interior and with
boundary of Haar measure $0$ such that
\begin{equation}
   \gL \, = \, \oplam(W) \, := \,
            \{ \pi^{}_G(y) \mid y\in\tilde{M},\, \pi^{}_J (y) \in W\}
       \, = \,  \{x \in M \mid x^* \in W \} \, .
\end{equation}

The word {\em regular\/} refers to the assumption on $W$
that its boundary, $\partial W$, has measure $0$. The assumption
that $W$ has non-empty interior is the non-trivial case of what
we are about to consider, and assures that the set
$\gL$ is relatively dense in $G$ (though this is not necessary for
what follows).

A {\em weighted subset\/} of $G$ is a subset $\gL$ of $G$
together with a (bounded) mapping $w \! :\gL \longrightarrow \CC$,
which we will assume to be bounded throughout this paper.

\begin{theorem} \label{cp-thm}
Let $(G,J,\tilde{M})$ be a cut and project scheme with
LCA groups $G$ and $J$, where $G$ is also $\sigma$-compact.
Let $W\subset J$ be relatively compact with boundary of Haar measure 
$0$.
Let $f$ be any complex-valued function on $J$ which is supported
and continuous on $\overset{\_\!\_\!\_\!\_}{W}$. Define a weighting $w$ 
on
$\oplam(W)$ by $w(x) := f(x^*)$. Let $\omega$ be the Dirac
comb\/ $\omega = \sum_{x\in \gL} w(x) \delta_x$.
Then, for any van Hove sequence $\{A_n\}$ in $G$,
the corresponding diffraction measure $\hat{\gamma}^{}_\omega$
is translation bounded, positive and pure point. In particular,
the weighted model set $(\oplam(W),w)$ is pure point diffractive.
\end{theorem}

This result, in various degrees of generality, has been proved by A.~Hof
\cite{Hof}, B.~Solomyak \cite{Boris}, and M.~Schlottmann \cite{Martin2}.
These proofs are all based upon the pointwise ergodic theorem for
uniquely ergodic dynamical systems, following or extending
Dworkin's argument \cite{Dworkin}.
Below, we present an alternative derivation which rests on our above
results. To our knowledge, this is the first proof
that links Theorem \ref{cp-thm} directly back to almost periodicity
rather than dynamical systems and ergodic theory. For more on
van Hove sequences, see the Appendix.

The strategy of the proof is to show that the hypotheses of Theorem
\ref{main-thm} are satisfied. For this, the main ingredient
is a version of Weyl's theorem on uniform projection.
In the form stated here, it is \cite[Thm.\ 6.2]{BM}, and the explicit 
proof
given there is based on the version of Weyl's theorem in \cite{Martin1}.
A new direct proof and a variant without the assumption
$\theta^{}_J (\partial W) = 0$ is derived in \cite{Moody2001}.

\begin{theorem} {\rm (Weyl)} \label{Weyl-thm}
Assuming the notation and prerequisites of Theorem $\ref{cp-thm}$,
\begin{equation}
    \lim_{n\to \infty}\, \frac{1}{\theta_G (A_n)}
    \sum_{x\in \mbox{\small $\curlywedge$}(W) \cap A_n} f(x^*)
    \; = \; C \int_W f \, {\rm d}\theta_J  \,,
\end{equation}
where $C$ is a positive contstant that depends only
on the choice of the Haar measures on $G$ and $J$ and
on the lattice $\tilde{M}$.  \qed
\end{theorem}

{\sc Proof of Theorem \ref{cp-thm}}:
We fix, once and for all, a van Hove sequence $\{A_n\}$ in $G$.
It is clear that $\omega$ is a translation bounded measure
since the set $\gL$ is uniformly discrete and the values
of $f$ are bounded in absolute value. Furthermore,
\begin{equation} \label{eq15}
     \eta^{}_n(z) \; = \;
     \frac{\big(\omega_n * \tilde{\omega}_n \big)(\{z\})}{\theta_G(A_n)}
     \; = \; \frac{1}{\theta_G(A_n)}
     \sum_{ \substack{x,y \in \gL \cap A_n \\ x-y = z} }
    {f(x^*)\overline{f(y^*)}} \, .
\end{equation}
Note that $z^*$ is well-defined since $z$ lies in $M$. Using the
function $f \,\overline{T_{z^*}(f)}$ in Weyl's theorem (where $T_a$
denotes the shift by $a$, so that $T_a(f)(x) = f(x-a)$) we obtain
\begin{equation}
    \eta(z) \; = \; \lim_{n \to \infty}
   \frac{\big(\omega_n * \tilde{\omega}_n \big) (\{z\})}{\theta_G(A_n)}
   \; = \; C \int_W f(u) \overline{f(u-z^*)} \, {\rm d}\theta_J(u)
\end{equation}
where the van Hove property of $\{A_n\}$ has been used to rewrite the
sum in (\ref{eq15}) so that Weyl's theorem can be applied.
In particular, the autocorrelation exists. Also,
since $\gL - \gL \subset \oplam(W - W)$ and $W-W$ has
compact closure, we obtain that $\gL - \gL$ is uniformly discrete,
and also the essential support of the autocorrelation
of $\omega$ is then uniformly discrete.

It remains to prove the $\ve$-almost periodicity. For this, it
suffices to show that, for all $\ve >0$, the set of points $z\in M$
for which $|\eta(0) - \eta(z)| < \ve$ is relatively dense.
Choose $\ve >0$. Then
\begin{eqnarray*}
    |\eta(0) -\eta(z)|  & = &
    \left|\, C \int_W
    f(u) \overline{\big( f(u) - f(u-z^*) \big)} \, {\rm d} \theta_J(u) 
\,\right|  \\
    & \le & C \int_W
    |f(u)| \,|f(u) - f(u-z^*)| \, {\rm d} \theta_J(u)  \, .
\end{eqnarray*}

Write $W = (W \cap (z^* +W)) \cup (W \backslash (z^* +W))$ and split
the integral accordingly. Using the
uniform continuity of $f$ on $W \cap (z^* +W)$, we can find
a neighbourhood $V_1$ of $0$ so that, for all $z^* \in V_1$,
the term with the integral over $W \cap (z^* +W)$ is bounded by $\ve/2$.
Using the boundedness of $f$ and the uniform continuity of
\begin{equation}
    \int_{W\backslash (u+W)} 1 \, {\rm d} \theta_J
    \; = \; \theta_J(W) -
    \big(\boldsymbol{1}_W*\tilde{\boldsymbol{1}}_W\big)(u)
\end{equation}
as a function of $u$,
we can find a second neighbourhood $V_2$ of $0$ on which the term
with the integral over $W\backslash (z^*+W)$ is  bounded by
$\ve/2$. Thus we obtain $|\eta(0) - \eta(z)| < \ve$
for all $z \in \oplam(V_1 \cap V_2)$. Since this set is relatively
dense ($V_1 \cap V_2$ is open) and a subset of $P_\ve$, the
latter is also relatively dense. \qed

\smallskip
Let us now come back to the discussion of Section \ref{sec2}, more 
precisely
to the structure of the sets $P_{\ve}$ of $\ve$-almost periods. If 
$0<\ve<1$,
the set $W=\overline{\varphi(P_{\ve})}$ has non-empty interior 
according to
Lemma \ref{model1}. This now has the immediate consequence that each 
such
$P_{\ve}$ is a subset of a model set, i.e., $P_{\ve} \subset \oplam(W)$ in 
the
setting of Eq.~(\ref{model-def}). Since we also have
$$ \overline{\varphi(P_{\ve})} \; \subset \;
     \{t\in H\mid \bar{\varrho}(t,0)\le \ve\} \, , $$
where $\bar{\varrho}$ means the (pseudo-)metric on $H$ induced
by $\varrho$,
the difference between $\oplam(W)$ and $P_{\ve}$ is a
subset of $\{t\in L\mid \varrho(t,0)=\ve\}$.
If our original set $S$ itself is a model set, then this property of the
$P_{\ve}$ relates to the fact that also $S-S$ is a model set 
\cite{Moody},
which, in turn, contains the $P_{\ve}$ for all $0<\ve<1$.
Moreover, since $L$ is countable, a difference between $\oplam(W)$ and 
$P_{\ve}$
is possible for at most countably many values of $\ve$. For all other 
values,
the sets $P_{\ve}$ {\em are\/} the model sets determined by their 
closures,
so that $\gD^{\rm ess}$ is actually the union of an ascending sequence
of model sets.

\section{Diffraction from visible lattice points}

Let $\gG$ be a lattice in $\RR^n$. If a basis is chosen, we can
write the points of $\gG$ as $t=(t_1, \dots ,t_n)$ with $t_i\in\ZZ$,
abbreviating the corresponding linear combination. Such a point is
called {\em visible\/} if $\gcd(t_1, \dots ,t_n) = 1$. The property
of visibility does not depend on the lattice basis chosen. The set
$V=V(\gG)$ of visible points of $\gG$ is uniformly discrete in $\RR^n$,
but it is not relatively dense. In fact, as follows from
\cite[Prop.~4]{BMP}, relative denseness cannot
be restored by adding points of density 0.
Nonetheless, if we fix a natural sequence $\cA$, e.g., an increasing
sequence of balls around the origin of $\RR^n$, we have:
\begin{prop} \label{vispo}
    The set\/ $V(\gG)$ is pure point diffractive, i.e., the Fourier
    transform of the autocorrelation of the Dirac comb\/
    $\omega^{}_V := \sum_{x\in V(\gG)} \delta_x$ is a positive
    pure point measure.
\end{prop}
Proposition \ref{vispo} was originally established in \cite{BMP}. Here, 
we offer an alternative proof based on the fact that Axioms \ref{ax1},
\ref{ax2}, 3$^+$, and \ref{ax4}
are satisfied by $\omega^{}_V$. The visible points thus show that the
regime of our axioms goes considerably beyond ordinary model sets
(which are always relatively dense).

We need the fact, established in \cite[Thms.\ 1 and 2]{BMP}, that the
autocorrelation of $\omega^{}_V$, for the sequence $\cA$ of increasing
balls around $0$,  exists and is a pure point measure
$\gamma^{}_{\omega^{}_V}$ supported on $\gG$. Its coefficient at
$t\in\gG$ is
\begin{equation} \label{vis-auto}
   \eta(t) \; = \; \gamma^{}_{\omega^{}_V} \big( \{t\} \big) \; = \;
   {\rm dens}(\gG)\, \xi(n) \prod_{p \, | \, {\rm cont}(t)}
   \Big( 1 + \frac{1}{p^n - 2} \Big)
\end{equation}
where ${\rm dens}(\gG)$ is the density of $\gG$,
${\rm cont}(t) := \max \{k\in\ZZ \mid t \in k \gG\}$, and
$$ \xi(n) \; = \; \prod_{p \; {\rm prime}} \Big( 1 - \frac{2}{p^n} 
\Big). $$
We also note that
$$ {\rm dens}(V) \; = \; \eta(0) \; = \;
    {\rm dens}(\gG)\, \xi(n) \prod_{p\; {\rm prime}}
    \Big( 1 + \frac{1}{p^n - 2} \Big)
    \; = \; \frac{{\rm dens}(\gG)}{\zeta(n)} $$
where $\zeta(n)$ is Riemann's zeta function \cite[Prop.\ 6]{BMP}.

In the notation of Section \ref{sec1}, we obtain
$\gD^{\rm ess} = \gD = V - V = \gG$, and Axioms \ref{ax1},
\ref{ax2}, and 3$^+$
clearly hold. To establish Axiom \ref{ax4}, choose any $\ve > 0$. Then,
\begin{eqnarray*}
   t \in P_{\ve} & \iff &
      1 - \frac{\eta(t)}{\eta(0)} \, < \, \ve^2 \\
   & \iff & \eta(t) > \eta(0) (1-\ve^2) = {\rm dens}(V) (1-\ve^2).
\end{eqnarray*}
Since $\zeta(n) \xi(n) \prod_{p\; {\rm prime}} 
\big(1+\frac{1}{p^n-2}\big)=1$,
we may choose $N$ so that
$$ \zeta(n) \, \xi(n) \prod_{\substack{p \; {\rm prime} \\ p\le N}}
    \Big(1+\frac{1}{p^n-2}\Big) \; > \; 1 - \ve^2 \, . $$
If $p^{}_1,\dots,p^{}_k$ are all the primes $\le N$, then it is
immediate from (\ref{vis-auto}) that, for all
$t\in p^{}_1\cdot\ldots \cdot p^{}_k \gG$, we have
$$ \eta(t) \; = \; {\rm dens}(V) \, \zeta(n) \, \xi(n)
    \prod_{p \, | \, {\rm cont}(t)} \Big( 1 + \frac{1}{p^n - 2} \Big)
    \; > \;  {\rm dens}(V) \, (1-\ve^2) $$
so that $P_{\ve} \supset p^{}_1\cdot\ldots \cdot p^{}_k \gG$ is
relatively dense.     \qed

\smallskip
The set $F(k)$ of integers in $\RR$ which are free of $k$th power
factors for some fixed $k$ are also shown in \cite{BMP} to form a
pure point diffractive set. We can again see that this situation
is covered by our present setup. This time, the key piece of
information is that the autocorrelation for $F(k)$ is  supported
on $\ZZ$, and the coefficients are given by
\cite[Thm. 4]{BMP}
\begin{equation} \label{power-free}
    \gamma^{}_{\omega^{}_{F}} (t) \; = \;
    \xi(k) \prod_{p^k \, | \, t} \Big( 1 + \frac{1}{p^k - 2} \Big).
\end{equation}
Otherwise, the argument is similar.

Let us mention that neither of the two examples of this Section is a 
model set
as defined in Section \ref{sec4}, even though both can be described in a
cut and project scheme with the adele ring over $\ZZ$ as internal space.
The special situation met here is that, in each case, the corresponding 
window $W$
is compact, but happens to be the boundary of another set, i.e.,
$W=\partial W'$. As such, $W$ has empty interior, which lines up
with the set of visible points (or the set of $k$-free integers) not 
being
relatively dense. To make matters worse, $W$ has positive measure, see 
\cite{BMP}
for details. The concept of {\em weak model set\/} introduced in 
\cite{Moody2001}
includes such sets, but the original diffraction results of 
\cite{Hof,Martin2} do
not apply and the original proof of pure point diffraction in 
\cite{BMP} is
rather painful. It is thus somewhat astonishing to observe that, in our 
present
setting with almost periodicity, examples such as the visible points 
turn
out to be relatively harmless.

\section{Fibonacci chain as model set}

{}For many tilings and point sets of importance, an explicit cut and 
project
scheme is known, or even is the only known way to define it. However, 
perhaps
the most famous ones, such as the Penrose tiling, are defined through
substitution rules or through perfect matching rules, and the additional
possibility to describe it by projection is still somewhat mysterious.
In particular, it is not clear what the intrinsic properties are which
lead to the projection picture.
Our above approach, at least in principle, allows to determine the 
internal
group and thus the most natural cut and project scheme for a given 
pattern.
We will explain this for the one-dimensional example of the Fibonacci 
chain,
which already shows the key features. Strictly speaking, we will make 
use
of the known cut and project scheme for it and only use the AC topology
to justify the choice of internal space {\em a posteriori}. It is 
possible
to avoid this, but then we would have to use the full machinery of
substitution systems to evaluate the autocorrelation, which is not 
really
enlightening.

The Fibonacci chain can be defined by means of the primitive two-letter
substitution rule $\sigma$ which sends $a\to ab$ and $b\to a$. A 
bi-infinite
sequence can easily be obtained as a fixed point of $\sigma^2$.
With $w^{}_1 = aa$ and $w^{}_{n+1} := \sigma^2 (w^{}_n)$, one obtains
the sequence
$$ a|a \; \stackrel{\sigma^2}{\longrightarrow} \;
    aba|aba \; \stackrel{\sigma^2}{\longrightarrow} \;
    \dots \; \xrightarrow{\, n\to\infty\,} \; w = \sigma^2(w)
$$
where $|$ denotes the origin or reference point chosen, and convergence
is in the obvious product topology on $\{a,b\}^\ZZ$.
Both letters occur with positive
frequency in the word $w$. A natural choice for a geometric
realization is an interval of length $\tau = \big( 1+\sqrt{5}\, \big)/2$
for $a$ and another of length 1 for $b$. Let us now consider the set
$\gL$ of left endpoints of this arrangement which is the ubiquitous
Fibonacci point set. Its difference set $\gD=\gL -\gL$ is a subset of
the ring $\ZZ[\tau]$ of integers in the quadratic field 
$\QQ\big(\sqrt{5}\,\big)$.
On the other hand, $\langle 1,\tau \rangle^{}_{\ZZ} = \ZZ[\tau]$, whence
we get $L=\ZZ[\tau]$ because $1$ and $\tau$ are certainly in 
$\gD_{}^{\rm ess}$.

It is known that $\gL$ is a model set \cite{Moody}. It can be described 
as
$$  \gL \; = \; \oplam(W) \; = \; \{ x \in \ZZ \mid x^* \in W \} $$
where $W = (-1,\tau-1]$ is a half-open interval and $*$ denotes 
algebraic
conjugation in $\QQ\big(\sqrt{5}\,\big)$
as defined by $\sqrt{5}\mapsto -\sqrt{5}$.
Also, $\gL-\gL$ is a model set, this time with window $[-\tau,\tau]$.
Using this and Weyl's theorem on
uniform distribution, see Theorem~\ref{Weyl-thm} above, one obtains a
closed formula for the autocorrelation of $\gL$. If $z\in \ZZ[\tau]$, 
one has
\begin{equation} \label{Fibo-corr}
    \eta(z) \; = \; \frac{{\rm dens}(\gL)}{{\rm vol}(W)}
         \int_{\RR}  1^{}_{W} (\xi)  1^{}_{W-z^*} (\xi) \, {\rm d}\xi \, 
,
\end{equation}
where $1^{}_W$ denotes the characteristic function of the set $W$,
and $\eta(z)=0$ otherwise. So, $\varrho(z,0)$ small means $\eta(z)$ 
close
to $\eta(0)$, and the latter is tantamount to saying that $z^*$ is small
in the usual topology of $\RR$. In other words, the Hausdorff completion
of $L=\ZZ[\tau]$ in the AC topology is $\RR$, and $\varphi$ is nothing
but the star map in this setting. This gives the promised a posteriori
justification for the standard cut and project setup used for the 
Fibonacci
chain.

It is perhaps remarkable that the star map is totally discontinuous in 
the
original topology of $G=\RR$, while it becomes uniformly continuous in 
the
new AC topology, and, with hindsight, this could have been a guide to 
finding
the extra topology much earlier.

Let us add that very much the same procedure works for many of the 
famous
standard examples, see \cite{MB} and references given there for the
zoo tamed so far. The above approach will always yield a setting that is
equivalent to what is called the minimal embedding case in previous 
works.
One particularly interesting case is the rhombic Penrose tiling, whose 
vertex
set is commonly described either in terms of a projection from four or 
from
5 dimensions. The former is the minimal embedding case, but one then 
has to
view the set as a multiple component model set. In our new setting, this
would be reflected by $H$ being the direct product $\RR^2 \times C_5$,
where $C_5$ is the cyclic group of order 5. If one insists on a 
Euclidean
internal space, one can embed $C_5$ into another copy of $\RR$ which
gives the common description with 3-space as internal group. However, 
this
is not a cut and project scheme in the strict sense because the 
projection
of the lattice into internal space is no longer dense, see \cite{Moody}
for details.

\section{Period doubling chain as model set}

Not all model sets in Euclidean space (i.e., with $G=\RR^n$) can be 
described
with an internal group $H$ that is also Euclidean, and there are 
relevant
examples \cite{BMS} where one has to go beyond. With our above setting, 
one
is able to actually determine the internal space also in this situation.
We will now illustrate this for the example of the period doubling 
point set,
where $H$ will turn out to be 2-adic. The resulting point set is an
example of a Toeplitz sequence and has been studied before, see 
\cite{HME}
and references therein. What is new here is the derivation of the 2-adic
integers from intrinsic data, extracted via the autocorrelation 
topology.

Let $\sigma$ be the primitive two-letter substitution defined by
$a \to ab$ and $b\to aa$. This is the well-known period doubling
substitution \cite[p.~301]{Allouche}. As above, a bi-infinite sequence
can be obtained as fixed point of $\sigma^2$,
$$  b|a \; \stackrel{\sigma^2}{\longrightarrow} \;
    abab|abaa \; \stackrel{\sigma^2}{\longrightarrow} \;
    \dots \; \xrightarrow{\, n\to\infty\,} \; w = \sigma^2(w)
$$
where $|$ denotes again the origin or reference point.
Note that this infinite word $w$ is
a so-called singular 2-cycle of $\sigma$ because $w$ and $\sigma(w)$
differ only at the position to the left of the marker, although they
are locally indistinguishable. This is, however, of no relevance to
our analysis because both $w$ and $\sigma(w)$ are repetitive.

Let us realize the symbolic sequence as a `coloured' point set in $\RR$
by using an interval of length 1 for both symbols, $a$ and $b$, but
different types (colours) of endpoints. Let
$\gL_a$ and $\gL_b$ denote the sets of left endpoints of all
$a$- and $b$-intervals, so that $\gL_a \cup \gL_b = \ZZ$.
Since $\sigma^2$ maps $a$ to $abaa$ and
$b$ to $abab$, we obtain the equations
\begin{eqnarray*}
   \gL_a & = & 4 \gL_a \cup (4 \gL_a + 2) \cup (4\gL_a + 3)
                 \cup 4\gL_b \cup (4\gL_b + 2) \\
   \gL_b & = & (4\gL_a+1) \cup (4\gL_b+1) \cup (4\gL_b + 3)\, .
\end{eqnarray*}
This can be simplified to the decoupled equations
\begin{eqnarray*}
   \gL_a & = & 2\ZZ \cup (4\gL_a + 3)  \\
   \gL_b & = & (4\ZZ + 1) \cup (4\gL_b + 3)
\end{eqnarray*}
which, by iteration, leads to the solution
\begin{eqnarray} \label{a-1}
   \gL_a & = & \bigcup_{n\ge 0} \big(2\cdot 4^n \ZZ + (4^n - 1) \big) \\
   \gL_b & = & \bigcup_{n\ge 1} \big(4^n \ZZ + (2\cdot 4^{n-1} - 1)
               \big) \, \cup \, \{ -1 \} , \label{b-1}
\end{eqnarray}
The r\^{o}le of $\{-1\}$ is exceptional since the other fixed point
of $\sigma^2$ (where the $b$ at $-1$ is replaced by an $a$, see above)
has $-1\in\gL_a$ rather than in $\gL_b$. Fortunately, this is of no
relevance to the autocorrelation, and we will simply suppress the
point $(-1)$ in the following calculations.
Since both point sets consist of disjoint unions of translates of
lattices whose lattice constants increase as powers of 4, they are
examples of what is usually called a limit periodic point set
\cite[p.~160]{Franz}.

We are ultimately interested in understanding $\gL_a$ and $\gL_b$ in
terms of our above concepts, which needs the knowledge of the
autocorrelation. We will sketch how this can be calculated in this case,
without going too much into detail. The idea is to approximate the
sets by periodic point sets obtained from {\em finite\/} unions in
(\ref{a-1}) and (\ref{b-1}), and to use Poisson's summation formula for
the approximants to determine their diffraction measures via the usual
intensity formula, i.e., by the fact that the diffraction intensities
of a periodic structure are the absolute square of the
corresponding amplitudes
(or Fourier-Bohr coefficients). The autocorrelation is then obtained
by inverse Fourier transform, and each step can easily be made rigorous
because all approximations involved here are appropriately convergent.

So, let $\omega_{a} = \sum_{x\in\gL_{a}} \delta_x$ be the Dirac
comb attached to $\gL_a$, and similarly for $\omega_{b}$, however
with the point measure at $x=-1$ removed.
In view of the terms in (\ref{a-1}) and (\ref{b-1}),
Poisson's summation formula gives
\begin{eqnarray*}
    \hat{\omega}_a & = &
    \sum_{n\ge 0} \frac{e^{-2\pi i k (4^n - 1)}}{2\cdot 4^n}
    \,\delta^{}_{\ZZ/2\cdot 4^n}  \\
    \hat{\omega}_b & = &
    \sum_{n\ge 1} \frac{e^{-2\pi i k (2\cdot 4^{n-1} - 1)}}{4^n}
    \,\delta^{}_{\ZZ/4^n}
\end{eqnarray*}
where $\delta^{}_{\gG}$ is the uniform Dirac comb of a lattice $\gG$
and $k$ is the wave number, i.e., the variable of the function in
front of the point measures. A short reflection shows that the point
measures are supported on the set
\begin{equation} \label{Bragg-support}
   F \; = \; \big\{\, \frac{m}{2^r} \mid
    (r=0 , \, m\in\ZZ ) \text{ or } (r\ge 1, \, m \text{ odd})\, \big\}
\end{equation}
which we will use as a parametrization of the support from now on.

{}For a superposition $\omega = \alpha\, \omega_a + \beta\, \omega_b$,
the diffraction would then be
\begin{equation} \label{pd-diffrac}
   \hat{\gamma}^{}_{\omega} \; = \;
   \sum_{k\in F} | \alpha A(k) + \beta B(k) |^2 \,
   \delta_k
\end{equation}
where the amplitudes, after a straightforward though somewhat lengthy
calculation, turn out to be
\begin{eqnarray}
   A(k) & = &  \frac{2}{3} \, \frac{(-1)^r}{2^r} \, e^{2\pi i k}
   \\ [2mm]
   B(k) & = & \delta^{}_{r,0} - A(k)
\end{eqnarray}
for any $k\in F$ of (\ref{Bragg-support}), and with $\delta^{}_{r,0}$
being Kronecker's function.
Choosing $\beta=0$ and $\alpha=1$ in (\ref{pd-diffrac}), one
obtains, by inverse Fourier transform, the autocorrelation of
$\omega_a$ as
$$  \gamma^{}_a \; = \; \sum_{z\in \ZZ} \eta^{}_a (z)\, \delta^{}_z $$
with the autocorrelation coefficients
\begin{equation} \label{a-coeff}
   \eta^{}_a(z) \; = \; \frac{2}{3} \cdot \begin{cases}
    \; \; 1 & \text{if $z=0$,} \\
    \Big( 1 - \frac{1}{2^{\ell + 1}} \Big) &
    \text{if $z=(2m+1) 2^{\ell}$, $\ell\ge 0$.}
    \end{cases}
\end{equation}
For completeness, let us add that the other autocorrelation
coefficients are
$$ \eta^{}_b (z) \; = \; \eta^{}_a (z) - \frac{1}{3} $$
for all $z\in\ZZ$, while those for the correlation between
the two types of points turn out to be
$$ \eta^{}_{ab} (z) \; = \; \eta^{}_{ba} (z) \; = \;
    \frac{2}{3} - \eta^{}_a (z)  $$
because we have $\eta^{}_a (z) + \eta^{}_b (z) +
\eta^{}_{ab} (z) + \eta^{}_{ba} (z) = 1 $ and the
general inversion symmetry follows from the palindromicity of
the period doubling sequence, see \cite{Baake-palin}.

Let us return to $\eta^{}_a (z)$ of (\ref{a-coeff}). If
$\ve > 0$ is given, then we have
$$ \big|\, \eta^{}_a(0)  - \eta^{}_a
    \big((2m+1)2^{\ell}\big) \, \big| \, < \, \ve $$
for all $\ell \ge \ell_0$ (with suitable $\ell_0$), and the
$\ve$-almost periods of $\gamma_a$ are thus relatively dense.
This shows, without using the model set description, that
the set $\gL_a$ (and also $\gL_b$) conforms to Axiom \ref{ax4},
and hence to all Axioms needed to apply Theorem \ref{main-thm}, so
the set is pure point diffractive. We have actually already seen
the explicit diffraction formula in (\ref{pd-diffrac}).

The translation invariant pseudo-metric (\ref{metric1}) on $\ZZ$
is given by
$$  \varrho^{}_a (z,0) \; = \;
     \left| 1 - \frac{\eta^{}_a(z)}{\eta^{}_a(0)}\right|^{1/2}
    \; = \; \begin{cases} \; 0 & \text{if $z=0$,} \\
            \; 2^{-(\ell+1)/2} &
    \text{if $z=(2m+1) 2^{\ell}$, $\ell\ge 0$.}
    \end{cases} $$
This clearly defines the 2-adic topology on $L = \ZZ$ and we
conclude that the internal group $H$ of Section~\ref{sec2}
is the 2-adic completion $\overline{\ZZ}_2$ of $\ZZ$.
The autocorrelation coefficient $\eta^{}_b (z)$ leads to the
same conclusion.

In this setting, $\overset{\_\!\_}{\gL_a}$ and 
$\overset{\_\!\_}{\gL_b}$ are
compact subsets of $\overline{\ZZ}_2$ which satisfy
\begin{eqnarray*}
    \overset{\_\!\_}{\gL_a} \cup \overset{\_\!\_}{\gL_b}  & = & 
\overline{\ZZ}_2 \\
    \overset{\_\!\_}{\gL_a} \cap \overset{\_\!\_}{\gL_b}  & = & \{ -1 \} 
.
\end{eqnarray*}
{}From this, we can obtain the model set description
$$  \gL_a \; = \; \oplam(W_a) \; , \quad
     \gL_b \; = \; \oplam(W_b)  $$
with $W_a = \overset{\_\!\_}{\gL_a} \setminus \{-1\}$ and
$W_b = \overset{\_\!\_}{\gL_b}$. If we had chosen the other fixed point
of $\sigma^2$, $\{-1\}$ would have gone to $W_a$ instead
of $W_b$.

\section{Extensions: Almost periodic measures}  \label{extensions}

Our restriction so far to Dirac combs $\omega$ with an autocorrelation 
that
is a pure point measure supported on a closed and discrete set (or on
a uniformly discrete set, if we invoke Axiom $3^+$) was motivated by
the possibility of a constructive, direct proof of the diffraction
results, and on the explicit construction of an internal space for
the corresponding cut and project scheme. In physical applications,
such Dirac combs are valid idealizations, but one would also like to
know the nature of the diffraction measure for more general
translation bounded measures.

The first (and rather obvious) extension concerns the case that $\omega$
has the form
\begin{equation} \label{conv-ext}
      \omega \; = \; g * \delta^{}_S
\end{equation}
where $\delta^{}_S$ is a Dirac comb of the original form considered in
Eq.~(\ref{comb1}) and $g$ is any finite measure, e.g., a (possibly 
continuous)
$L^1$-function. This would describe a situation where the measure 
$\omega$
emerges from adding up `local' contributions at the positions of the
set $S$. In this case, the autocorrelation exists whenever that for 
$\delta^{}_S$
does, and one obtains the convolution
$$  \gamma^{}_{\omega} \; = \; (g * \tilde{g}) * \gamma^{}_S \, . $$
The diffraction measure, by the convolution theorem, reads
$$  \hat{\gamma}^{}_{\omega} \; = \;
     | \hat{g} |^2 \cdot \hat{\gamma}^{}_S \, . $$
So, if a translation bounded measure $\omega$ admits a factorization
as in (\ref{conv-ext}), we can apply our previous analysis directly
to the discrete part of it, $\delta^{}_S$, and obtain the result on
the basis of the properties of $S$. So, if $\delta^{}_S$ satisfies
Axioms \ref{ax1}, \ref{ax2}, $3^+$ and \ref{ax4}, also the diffraction
measure $\hat{\gamma}^{}_{\omega}$ is pure point. In particular, this
situation applies to all lattice periodic measures which can always
be written as in (\ref{conv-ext}).

More interesting, however, is the question for a complete
characterization of translation bounded measures $\omega$ with
pure point diffraction measure. Such a characterization can be
given on the basis of \cite{GdeL}, however for the price of losing
the constructive part of our above analysis.

{}For this, we have to work with the space $\cM^{\infty}(G)$, but
equipped with other topologies. Let us first introduce the (local)
norm topology. Let $K$ be a compact neighbourhood of $0\in G$, which
we assume fixed from now on. We set
\begin{equation} \label{norm-top}
    \| \omega \|^{}_K \; := \; \sup_{x\in G} |\omega| (x+K)
\end{equation}
where $|\omega|$ denotes the total variation of $\omega$, which is
also translation bounded. It is clear that $\|.\|_K$ defines a norm,
and we call the corresponding topology the {\em local norm topology\/}. 
In this setting, a measure $\omega$ is translation bounded iff 
$\|\omega\|^{}_K < \infty$, see the remarks in \cite[p.\ 145]{Martin2}. 
It is also clear that any other compact neighbourhood, $K'$ say,
would lead to an equivalent norm, and hence to the same topology,
because $K'$ can be covered by finitely many translates of $K$
and vice versa. 

\smallskip \noindent {\sc Remark}
This norm makes $\cM^{\infty}(G)$ into a Banach space, compare 
\cite[Ex.\ 2.13 and Thm.\ 2.3]{GdeL}. 
To see this, one can check that the supremum norm in \cite{GdeL}, defined via
uniform partitions, is equivalent to our approach, and hence defines the
same topology. There are systematic reasons to prefer an approach via
uniform partitions or uniform coverings, if one also needs other norms of
$\ell_{}^{p}$-type. Since we only need the supremum norm (with $p=\infty$),
the simpler approach is sufficient. 

\smallskip
We say that a measure $\omega$ is {\em norm almost periodic\/} if, for 
all $\ve > 0$, the set
$\{\, t \mid \| \delta_t * \omega - \omega \|^{}_K < \ve \}$
of norm $\ve$-almost periods is relatively dense in $G$.
This type of definition for almost periodicity goes back to H.\ Bohr.
As was pointed out by Bochner \cite{Bochner},
almost periodicity is often equivalent to a certain relative
compactness criterion, which can prove more convenient. The 
equivalence of these two notions, in a Banach space say, requires the 
continuity of the group action in the corresponding topology. This is 
{\em not\/} the case in our situation here.

\begin{prop} \label{almost-per-1}
     Let\/ $\omega$ conform to Axioms\/ $\ref{ax1}$, $\ref{ax2}$, and\/ $3^+$.
     Then, its autocorrelation, $\gamma^{}_{\omega}$, satisfies
     Axiom~$\ref{ax4}$ if and only if it is norm almost periodic.
\end{prop}
{\sc Proof}: We have $\omega\in\cM^{\infty}(G)$ (Axiom \ref{ax1})
and $\gamma^{}_{\omega} = \sum_{z\in\gD_{}^{\rm ess}} \eta(z) \delta_z$
exists (Axioms \ref{ax2} and $3^+$), with $\gD_{}^{\rm ess}$ uniformly
discrete. There is nothing to be shown if
$\eta(0) = \gamma^{}_{\omega} (\{0\}) = 0$, so let us assume
$\eta(0) > 0$ from now on.
If $t\in L$, we have $t+L=L$. For arbitrary $x\in G$,
we then obtain

\begin{eqnarray*}
     | \delta_t * \gamma^{}_{\omega} - \gamma^{}_{\omega} | \, (x+K) & = 
&
      \Big|\sum_{z\in L} \big(\eta(z-t) - \eta(z)\big)\, \delta_z\Big|
      \, (x+K) \\
    &\le & \sum_{z\in L} \big| \eta(z-t) - \eta(z) \big|
       \, \delta_z (x+K) \\
    & = & \sum_{z\in F_x} \big| \eta(z-t) - \eta(z) \big|
\end{eqnarray*}
where $F_x = \big(\gD_{}^{\rm ess} \cup (t + \gD_{}^{\rm ess})\big)
\cap (x+K)$ is a {\em finite\/} set. Its cardinality is bounded above by
a constant $C$, uniformly in $x$ and $t$, due to the uniform
discreteness of the set $\gD_{}^{\rm ess}$. In fact, one sees from
this calculation that the inequality here is actually equality, something 
that we will use in the converse part of the proof below.

Assume that Axiom~\ref{ax4} is in force.
Fix $\ve>0$ and set $\ve' = \ve/\big(C\eta(0)\sqrt{2}\,\big)$. For any
$t\in P_{\ve'}$, we then obtain from Lemma~\ref{krein-bound} that
$|\eta(z\pm t) - \eta(z)| < \ve/C$ for all $z\in G$, where we have
used the symmetry of $P_{\ve'}$ (see Fact~\ref{inclusion}).
Consequently,
$$ \| \delta_t * \gamma^{}_{\omega} - \gamma^{}_{\omega} \|^{}_K
     \; < \; C \sup_{z\in L} \big| \eta(z-t) - \eta(z)\big|
     \; \le \; \ve \, . $$
But $P_{\ve'}$ is relatively dense according to Axiom \ref{ax4}, and
then so are the norm $\ve$-almost periods of $\gamma^{}_{\omega}$, for
all $\ve>0$, which proves one direction.

Conversely, let $\gamma^{}_{\omega}$ be norm almost periodic and
set
$$ Q_{\ve} \; := \; \{t\in L\mid \|\delta^{}_t*\gamma^{}_{\omega}
     - \gamma^{}_{\omega}\|^{}_K < \ve \} . $$
This is again a symmetric set
because $\|\delta_x * \mu\|^{}_K = \|\mu\|^{}_K$ for all
$\mu\in\cM^{\infty}(G)$ and $x\in G$.
Observe that, for all $z\in G$ and $t \in L$, we have
$$ | \eta(z \pm t) - \eta(z)| \; \le \;
     \| \delta^{}_{\mp t} * \gamma^{}_{\omega} -
        \gamma^{}_{\omega}\|^{}_K  \; = \;
      \| \delta^{}_{t} * \gamma^{}_{\omega} -
        \gamma^{}_{\omega}\|^{}_K  \, .$$
Fix $\ve>0$ and set $\ve' = \eta(0) \ve^2$ which is still positive.
For all $t\in Q_{\ve'}$, we get
$$  \varrho(t,0) \; = \;
     \left( \frac{|\eta(t) - \eta(0)|}{\eta(0)}\right)^{1/2}  \le \;
     \left( \frac{\| \delta^{}_{t} * \gamma^{}_{\omega} -
       \gamma^{}_{\omega}\|^{}_K}{\eta(0)}\right)^{1/2}
       < \;  \ve $$
The set $Q_{\ve'}$ is relatively dense which implies that also
$P_{\ve}$ is relatively dense. This being true for all $\ve>0$,
Axiom~\ref{ax4} is satisfied.   \qed

\smallskip
To be able to use the results of \cite{GdeL}, we have yet to introduce
another topology on $\cM^{\infty}(G)$, called the {\em product topology}
in \cite[Ex.\ 2.15]{GdeL}. If $C^{}_U(G)$ denotes the space of bounded,
uniformly continuous functions on $G$, a measure $\mu\in\cM^{\infty}(G)$
is identified with an element of the Cartesian product space
$\big[ C^{}_U(G)\big]^{\cK(G)}$ via
$$  \mu \; = \; \{ g * \mu \} ^{}_{g\in\cK(G)}  $$
where each $g*\mu$ is uniformly continuous and bounded by
$\| g*\mu \|^{}_{\infty} \le \|g\|^{}_{\infty}
   \sup_{t\in G} |\mu| \big(t + {\rm supp}(g)\big)$,
hence an element of $C^{}_U(G)$.
This way, we have $\cM^{\infty}(G)\subset\big[ C^{}_U(G)\big]^{\cK(G)}$
which makes $\cM^{\infty}(G)$ into a locally convex topological vector
space, in the relative topology. A fundamental system of semi-norms
is provided by $\|\mu\|_g := \|g*\mu\|_{\infty}$, with
$g\in \cK(G)$. Unfortunately, $\cM^{\infty}(G)$
is not a complete subspace of $\big[ C^{}_U(G)\big]^{\cK(G)}$, but it
is {\em bounded closed\/}, i.e., every closed and bounded subset of
$\cM^{\infty}(G)$ is also complete in $\big[ C^{}_U(G)\big]^{\cK(G)}$,
see \cite[Thm.\ 2.4 and Cor.\ 2.1]{GdeL}.
Let us add that, due to the product structure, a subset $M$ of
$\cM^{\infty}(G)$ is bounded iff, for all $g\in\cK(G)$,
we have $\|g*\mu\|^{}_{\infty}\le C$ for all $\mu\in M$,
with a constant $C=C(g)$ which does not depend on $\mu$, compare
\cite[Ex.\ II.4.7 c]{Bou}.

A measure $\omega\in\cM^{\infty}(G)$ is then called {\em strongly almost
periodic\/} if, for each open neighbourhood $V$ of $0$ in the product topology, 
$\{t \in G \,  | \,  \delta_t *\omega  - \omega \in V \}$ is relatively dense in $G$.
In fact, in the case of the product topology, this is equivalent to 
requiring the set $\{ \delta_t * \omega \mid t\in G\}$ to be relatively 
compact. The use of the word `strong' in this context comes from \cite{GdeL}.

\begin{lemma} \label{norm-implies-product}
    If\/ $\omega\in\cM^{\infty}(G)$ is a norm almost periodic measure,
    it is also strongly almost periodic.
\end{lemma}
{\sc Proof}: First, let $g\in\cK(G)$ and $\mu\in\cM^{\infty}(G)$. Then
we obtain, with $V:= {\rm supp}(g)$,
\begin{eqnarray*}
   \left| \int_G g(x-y) \, {\rm d}\mu(y) \, \right|   & \le &
   \int_G \big| g(x-y) \big| \, {\rm d} |\mu|(y)  \\
   & \le & \int_G \bs{1}^{}_V (x-y) \, \|g\|^{}_{\infty}
     \, {\rm d} |\mu|(y)  \\ [1mm]
   & = & \| g \|^{}_{\infty} \cdot |\mu| (x - V)  \\ [1mm]
   & \le &  \| g \|^{}_{\infty} \,
     \sup_{x\in G} |\mu| \big( x + (-V)\big) \\
   & \le & C_g \, \| g \|^{}_{\infty} \, \|\mu\|^{}_K
\end{eqnarray*}
where $C_g > 0$ is a constant which emerges from the observation that
the support of $g$ is compact and can thus be covered by finitely
many translates of $K$. So, we have
$$  \| g * \mu \|^{}_{\infty} \; = \;
     \sup_{x\in G} \left| \int_G g(x-y) \, {\rm d}\mu(y) \,\right|
     \; \le \; C_g \, \| g \|^{}_{\infty} \, \|\mu\|^{}_K \, . $$

Let us now choose an arbitrary finite family
${\mathcal F} = \{g^{}_1,\dots ,g^{}_n\}$ of non-zero functions
$g_i\in\cK(G)$ and an $\ve >0$. Then,
$$  V^{\mathcal F}_{\ve} (0) \; := \;
     \{ \mu\in\cM^{\infty}(G) \mid \sup_{1\le i \le n}
      \| g_i * \mu \|^{}_{\infty} < \ve \} $$
is a neighbourhood of the 0-measure in $\cM^{\infty}(G)$ in the
product topology. The neighbourhoods of this kind form a fundamental
system of neighbourhoods of the 0-measure.
Choose $\ve' = \ve / \max\{C_{g_i} \|g_i\|^{}_{\infty}\mid 1\le i\le 
n\}$.


If $t\in G$ is an $\ve'$-almost period of $\omega$ for the local norm topology then for $\mu:= \delta_t*\omega - \omega$
we have $||\mu||_K < \ve'$, whence $||g_i*\mu||_\infty <\ve$. It follows that $t$ is a
$V^{\mathcal F}_{\ve} (0)$-almost period of $\omega$ for the strong topology. The lemma
follows at once.
\qed

\smallskip
Let us now recall the following result from \cite{GdeL}, which
is proved on the basis of the
Bohr compactification of $G$, see \cite[Ch.\ 7]{GdeL}.
\begin{prop} \label{arga-thm}
   Let $\mu\in\cM^{\infty}(G)$ be a transformable measure,
   with associated measure $\hat{\mu}\in\cM^{\infty}(\hat{G})$.
   Then, $\hat{\mu}$ is a pure point measure if and only if
   $\mu$ is strongly almost periodic.
\end{prop}
{\sc Proof}: This follows from \cite[Thm.\ 11.2 and Cor.\ 11.1]{GdeL}
by an application of the inverse Fourier transform. \qed

\smallskip
Combining Proposition~\ref{arga-thm} with Fact~\ref{general-fourier},
we obtain:
\begin{theorem}  \label{final-thm}
   Let $\omega$ be a translation bounded complex measure on the
   $\sigma$-compact, locally compact Abelian group $G$ and assume
   that, w.r.t.\ an averaging sequence $\cA$,
   its autocorrelation measure $\gamma^{}_{\omega}$ exists.
   Then, $\gamma^{}_{\omega}$ is a translation bounded positive definite
   measure, whose Fourier transform, $\hat{\gamma}^{}_{\omega}$,
   exists and is a translation bounded positive measure on $\hat{G}$.

   Moreover, $\hat{\gamma}^{}_{\omega}$ is a pure point measure
   $($i.e., $\omega$ has pure point diffraction$)$ if and only
   if $\gamma^{}_{\omega}$ is strongly almost periodic.  \qed
\end{theorem}

Consequently, together with Proposition~\ref{almost-per-1} and
Lemma~\ref{norm-implies-product}, we get an independent derivation
of our main Theorem~\ref{main-thm}.

It is not clear to what extent norm and strong almost periodicity
coincide, as the product topology is weaker than the norm topology.
The demarkation apparently is connected with the difference between
Axioms~\ref{ax3} and $3^+$. It is possible though to come to a partial
converse of Lemma~\ref{norm-implies-product}. For that, we first
need a technical result.
\begin{lemma} \label{tech}
    Let $u,v$ be two arbitrary complex numbers which satisfy the
    simultaneous inequalities $ | u - b v | < \ve $ and\/
    $|v - a u | < \ve $, with some real numbers $0\le a,b \le 1$.
    Then, they also satisfy\/ $|u-v| < 3\,\ve$.
\end{lemma}
{\sc Proof}:
Note first that $1-ab\ge 0$, so that we have
\[
   \lvert v \rvert (1-ab) \; = \;
   \lvert v - abv \rvert \; \le \;
   \lvert v - au \rvert + a \lvert u - bv \rvert
   \; < \; 2\ve .
\]
Since $0\le 1-b \le 1 - ab$, we now obtain
\[
   \lvert u-v \rvert \; \le \; \lvert u - bv \rvert +
   \lvert bv - v \rvert \; < \; \ve + \lvert v \rvert
   (1-b) \; < \; 3\ve,
\]
which was to be shown. \qed

\smallskip
Before we state the partial converse, let us make one further 
observation, making use of the fact that the two possible definitions
of almost periodicity (via $\ve$-periods and via relative compactness)
are equivalent in the product topology.
Strong almost periodicity of a measure $\mu\in\cM^{\infty}(G)$ means
relative compactness of the set $M_{\mu}=\{\delta_t*\mu\mid t\in G\}$ in
the product topology, which coincides with precompactness here. So, for
every finite family ${\mathcal F}$ of functions and for every $\ve>0$,
there are {\em finitely\/} many measures $\{\nu^{}_1,\dots,\nu^{}_n\}$
such that
\begin{equation} \label{condition-measures}
    M_{\mu} \; \subset \;\, \bigcup_{i=1}^{n} \,
               V^{\mathcal F}_{\ve}(\nu^{}_i)
\end{equation}
with the neighbourhoods $V^{\mathcal F}_{\ve}$ defined in analogy to 
above, i.e., they are translates of those around $0$.
Without loss of generality, we may assume that all the $\nu^{}_i$ are
actually translates of $\mu$, i.e., of the form
$\nu^{}_i = \delta_{t_i} * \mu$ for some $t_i\in G$.
Observing now that
$$ V^{\mathcal F}_{\ve}(\delta_{t}*\mu) \; = \;
    \delta_t * V^{\mathcal F}_{\ve}(\mu)\, , $$
we can rephrase the relative compactness condition and
(\ref{condition-measures}) by saying that there are finitely
many translations $\{t_1,\dots,t_n\}$ in $G$ such that
$$ M_{\mu} \; \subset \;\, \bigcup_{i=1}^{n} \,
      \big(\delta_{t_i} * V^{\mathcal F}_{\ve}(\mu)\big) . $$
This formulation is of advantage for the next result, and can also
be used to see that strong almost periodicity of $\mu$ implies that
all functions $g*\mu \in C^{}_U(G)$, with $g\in\cK(G)$, are
(uniformly) almost periodic functions on $G$, see \cite[p.\ 20]{Zaidman}
for a definition.

\begin{prop}  \label{partial-converse}
    Let $\mu\in\cM^{\infty}(G)$ be a pure point measure that is
   supported on a set $\gD$ for which $\gD - \gD$ is uniformly discrete.
   If $\mu$ is strongly almost periodic, then it is also norm almost 
periodic.
\end{prop}
{\sc Proof}: We can assume $\mu = \sum_{z\in\gD} \eta(z)\delta_z$
with $\eta(z)\neq 0$ for all $z\in\gD$, i.e., we assume
$\gD = \gD_{}^{\rm ess}$.
Since $\gD-  \gD$ is uniformly discrete, we can choose a compact 
neighbourhood $S$ of $0\in G$ such that, for all $u \in G$, each of 
$\gD \cap (u + S)$ and $(\gD -\gD) \cap (u + (S-S))$ contains at most one 
point. Without loss of generality, we may also assume that $S$ is
symmetric, i.e., $S=-S$, which simplifies some of the calculations.

Let us fix a function $g\in\cK(G)$ with values in the interval $[0,1]$,
support $S$, and $g(0)=1$. We can now apply our above
reasoning on strong almost periodicity to neighbourhoods of type
$V^{\{g\}}_{\ve /3}$, i.e., for all $\ve>0$, there exist
{\em finitely\/} many translations $\{t_1,\dots,t_n\}$ such that,
for all $t\in G$, we have
\begin{equation} \label{strong-ap}
    \| g * (\mu - \delta_{t-t_i} * \mu) \|^{}_{\infty} \; < \; \ve /3
\end{equation}
for (at least) one of the $t_i$. Here, the supremum is taken over the 
function
\begin{eqnarray*}
    \lefteqn{\big|\big(g*(\mu -\delta_{t-t_i}*\mu)\big)(x)\big|}\\[1mm]
    & = & \Big| \sum_{z\in\gD}\eta(z) g(x-z) \, -
          \sum_{z\in\gD}\eta(z) g(x-z-(t-t_i))\, \Big| \\[1mm]
    & = & \Big| \sum_{z\in\gD\cap(x+S)}\eta(z) g(x-z) \; - \!\!\!
          \sum_{z\in\gD\cap(x-(t-t_i)+S)}\!\!\!
          \eta(z) g(x-z-(t-t_i))\, \Big|
\end{eqnarray*}
where we made use of the symmetry of $S$ in the last step.
The choice of $S$ means that for each value of $x$,
there is at most one term from each of the sums. As $x$ varies over $G$,
one glides a $g$-shaped `dumb-bell' over $G$ with centres at
distance $t-t_i$. 
This distance is locally constant and changes only if
another $t_i$ has to be chosen.

If only one term is caught, $\eta(z)$ in the first sum say, then we may 
choose
$x=z$ to see that $|\eta(z)|<\ve/3$.
So, let us assume that $z^{}_1,z^{}_2$ are the two points of $\gD$
which are covered by the dumb-bell. Choosing first $x=z^{}_1$ and then
$x=z^{}_2$,
we obtain the simultaneous inequalities
$$ \big| \eta(z^{}_1) - b\, \eta(z^{}_2) \big| \; < \; \ve /3 \; , \quad
\big| a\, \eta(z^{}_1) - \eta(z^{}_2) \big| \; < \; \ve /3$$
where $a,b$ are the corresponding values of the function $g$, hence
numbers in the unit interval. Invoking Lemma~\ref{tech}, we get
$\big|\eta(z^{}_1) - \eta(z^{}_2)\big| < \ve$.

{}From this argument, we conclude that, for $\{t_1,\dots, t_m\}$ as
chosen for (\ref{strong-ap})
and for all $t \in G$, there is a $t_i \in \{t_1,\dots, t_m\}$ for which
\begin{equation}  \label{eq-32}
    \Big| \sum_{z\in\gD\cap(x+S)}\eta(z) \; - \!\!\!
    \sum_{z\in\gD\cap(x-(t-t_i)+S)}\!\!\! \eta(z) \, \Big| \; < \; \ve
\end{equation}
for all $x \in G$. Furthermore, at most one term from each sum is 
actually present.

We have to show that $\mu$ is norm almost periodic. Since
$\|.\|^{}_K$ and $\|.\|^{}_S$ are equivalent norms, they define
the same topology. So, we have to show that for each $\ve  >0$, the set
of $\ve$-almost periods with respect to $\|.\|^{}_S$ is relatively 
dense.

Take any $t \in G$ and choose the appropriate $t_i$ on the basis of 
\eqref{eq-32}. Consider
\[
  \mu - \delta_{t-t_i} * \mu \; = \;
  \sum_{z\in\gD} \eta(z) \delta_z \; -
    \sum_{z\in\gD} \eta(z)\delta_{t-t_i+z} \, .
\]
Let $x\in G$.  Restricting to $x+S$, we have
\[ 
  \sum_{z\in\gD\cap(x+S)} \eta(z) \delta_z\; - \!\!\!
    \sum_{z\in\gD\cap(x-(t-t_i)+S)}\!\!\! \eta(z) \delta_{t-t_i+z} 
\]
which reduces to
\[
   \eta(z_1) \delta_{z_1} - \eta(z_2)\delta_{t-t_i +z_2} \, ,
\]
if both $S$-sets actually meet $\gD$, and to just one or neither of the 
terms otherwise. Furthermore, $|\eta(z_1) - \eta(z_2)| < \ve$ if both terms 
are present, and similarly if just one is present.

Next, $z_2\, - z_1\in (\gD -\gD) \cap (t_i-t) + (S-S))$, which by the 
choice of $S$ has at most one point $p(t)$, which is independent of $x, z_1 
\,\textrm{and} \, z_2$. Thus, $z_2 = z_1 +p(t)$ and {\em for all\/} $z_1 \in \Delta$,
\[
  |\eta(z_1 + p(t)) - \eta(z_1)| \; < \; \ve \, ,
\]
where the first term is not there if $z_2$ does not exist.

Finally,
\begin{eqnarray}
   ||\delta_{-p(t)}*\mu \, - \mu  ||_S &=& ||\sum_{z\in\Delta} 
           (\eta(z)\delta_{-p(t) +z}\, -
    \eta(z))\delta_z ||_S \\
   &=& ||\sum_{z\in\Delta} (\eta(z+p(t)) -\eta(z))  \delta_z ||_S 
   \; < \; \ve \,  .
\end{eqnarray}

Thus, $-p(t)$ is an $\ve$-almost period and  $ t\in-p(t) + t_i + S-S$ 
shows that, with the compact set $K := \{t_1, \dots,  t_n\} + (S- S)$ 
and with the $\ve$-almost periods $AP_\ve$, $t \in AP_\ve + K$. 
Since $t$ was arbitrary, the $\ve$-almost periods are relatively 
dense.   \qed



\smallskip
Finally, to wrap up the results of this paper, we formulate:
\begin{theorem} \label{wrap-up}
    Let $\omega$ be a translation bounded measure on the $\sigma$-compact
    LCA group $G$ whose autocorrelation $\gamma^{}_{\omega}$ exists
    $($relative to an averaging sequence $\cA\,)$ and is a pure point
    measure with a support $\gD$ such that $\gD-\gD$ is still
    uniformly discrete $($in particular, $\omega$
    conforms to Axioms~$\ref{ax1}$, $\ref{ax2}$, and $3^+ \,)$. Then the
    corresponding diffraction measure $\hat{\gamma}^{}_{\omega}$ exists
    and the following statements are equivalent.
\begin{itemize}
    \item[{\rm (a)}] Axiom~$\ref{ax4}$ holds, i.e., the set $P_{\ve}$ 
of\/
              $(\ref{almost-periods})$ is relatively dense for all 
$\ve>0$.
    \item[{\rm (b)}] The autocorrelation $\gamma^{}_{\omega}$
                     is norm almost periodic.
    \item[{\rm (c)}] The autocorrelation $\gamma^{}_{\omega}$
                     is strongly almost periodic.
    \item[{\rm (d)}] The diffraction measure $\hat{\gamma}^{}_{\omega}$
                     is pure point.
\end{itemize}
\end{theorem}
{\sc Proof}: We have (a) $\Longleftrightarrow$ (b) by
Proposition~\ref{almost-per-1}, while
(c) $\Longleftrightarrow$ (d) follows from Proposition~\ref{arga-thm}.
Next, (b) $\Rightarrow$ (c) results from 
Lemma~\ref{norm-implies-product},
and the converse direction is Proposition~\ref{partial-converse}. \qed

\section*{Appendix: Averages and limits}

So far, we have fixed one averaging sequence $\cA$ out of a huge class 
of
possible sequences. This is not entirely satisfactory
because the resulting diffraction measure can, and generally will,
depend on this choice. Let us thus indicate how to restrict the choice 
of
$\cA$ in order to get rid of this problem, at least to some extent. The
method of choice is the use of so-called {\em van Hove sequences\/} 
which
were introduced long ago in statistical mechanics to deal with a
similar problem \cite[Ch.\ 2.1]{Ruelle}.

The key idea here is the following. If the Dirac comb (or measure)
$\omega$ is sufficiently `homogeneous' and if the sets of the sequence
$\cA$ grow in such a way that their boundaries are
suitably negligible (in measure) when compared
to the full sets, one should expect independence of the
autocorrelation of the actual choice of $\cA$. This is certainly the
case for lattices and general periodic structures in Euclidean spaces,
where a sequence of balls (with arbitrary centres) is as good as
one with cubes or similar shapes. To formalize this, consider two
{\em compact\/} subsets $B,K$ of $G$ and define the $K$-boundary of $B$ 
as
\begin{equation} \label{boundary}
   \partial^K B \; := \; \big((B+K)\setminus \overset{\,\circ}{B}\,\big)
   \,  \cup  \,
   \big((\,\overline{G\setminus B} - K)\cap B \, \big).
\end{equation}
This can be understood as some sort of $K$-thickening of
$\partial B$. In particular, $\partial^{\{0\}} B = \partial B$.
With this, $\cA$ of (\ref{sequence}) is called a (nested)
{\em van Hove sequence\/} if, for every compact $K\subset G$,
\begin{equation} \label{van-hove}
    \lim_{n\to\infty} \frac{\theta^{}_G
    (\partial^K \overset{\_\!\_}{A_n})}
    {\theta^{}_G (A_n)}  \; = \; 0\, .
\end{equation}
The meaning of (\ref{van-hove}) is that the boundary/bulk ratio
goes to zero (in measure) as $n\to\infty$, see \cite[Sec.\ 1]{Martin2}
for further details and consequences of this concept in the setting
of translation bounded measures. Note that general van Hove
sequences need not be nested, i.e., need not satisfy
$\overset{\_\!\_}{A_n} \subset A_{n+1}$, but, since $G$ is 
$\sigma$-compact, a general van Hove sequence always has a subsequence 
of the form $\cA' = \big(a_n + A_n\big)_{n\in\NN_0}$, for suitable 
$a_n \in G$, where the $A_n$ are properly nested and hence constitute
an averaging sequence $\cA$ according to \eqref{sequence}.

\begin{prop}
   Let $\gL$ be a regular model set in the sense of Section~$\ref{sec4}$,
   and let $\omega = \delta^{}_{\gL}$ be its uniform Dirac comb.
   Then, for any general van Hove sequence $\cA'$, the corresponding
   autocorrelation measure $\gamma^{}_{\omega}$ exists.
   Furthermore, it is independent of $\cA'$.
\end{prop}
{\sc Proof}: The result is true for point sets that lead to uniquely
ergodic dynamical systems under the action of $G$, see
\cite[Thm.\ 3.4]{Martin2}. A regular model set is such a point set
according to \cite[Thm.\ 4.5 (1)]{Martin2}. \qed

\smallskip
This result contains the cases of lattices and periodic point sets.
It can be extended to weighted model sets as we introduced them
above, but we will not expand upon this here.

Although the stated independence of the van Hove sequence chosen
is very plausible (or at least desirable) for many if not
most physical applications, the situation is not always as nice as
with regular model sets. Already in the case of the visible lattice
points, the above statement is clearly false because of the existence
of holes of arbitrary size. A sequence of balls centred around
the origin would give a different answer than one that moves around
chasing the holes. In other words, a certain degree of uniformity of
the limit is lost. Nevertheless, some partial result remains true.
If one restricts again to the {\em nested\/} form $\cA$ of the van
Hove sequence, existence of the autocorrelation and independence
of $\cA$ is regained, see \cite{BMP} for details.

\section*{Acknowledgements}

It is a pleasure to thank Tilmann Gneiting for a number of very
useful hints on the theory of positive definite functions and
their applications, and Aernout C.~D.~van Enter, Daniel Lenz
and Nicolae Strungaru for suggesting a number of improvements.
This work was supported by the German Research Council (DFG) and the
Natural Sciences and Engineering Research Council of Canada (NSERC),
and by the Volkswagen Foundation through the RiP program at Oberwolfach.

\end{document}